\documentclass[12pt,a4paper,reqno]{amsart}
\usepackage[headings]{fullpage}
\usepackage{amsthm}
\usepackage{amssymb}
\usepackage[nobysame]{amsrefs}
\usepackage[all]{xy}
\UseComputerModernTips
\numberwithin{equation}{section}
\theoremstyle{plain}
\newtheorem{theorem}{Theorem}
\newtheorem{proposition}[theorem]{Proposition}
\newtheorem{lemma}[theorem]{Lemma}
\newtheorem{corollary}[theorem]{Corollary}
\theoremstyle{definition}
\newtheorem{definition}[theorem]{Definition}

\theoremstyle{remark}
\newtheorem{remark}[theorem]{Remark}
\newtheorem*{acknowledgements}{Acknowledgements}
\newtheorem*{notations}{Notations and conventions}
\numberwithin{theorem}{section}
\newcommand{\F}{{\mathbb F}}

\newcommand{\Z}{{\mathbb Z}}
\newcommand{\Q}{{\mathbb Q}}
\newcommand{\C}{{\mathbb C}}

\newcommand{\thh}{{T\kern -1pt H\kern -1pt H}}
\newcommand{\tc}{{TC}}
\newcommand{\ti}{\tilde}
\newcommand{\wti}{\widetilde}
\newcommand{\lr}{\longrightarrow}
\newcommand{\cok}{\operatorname{cok}}
\newcommand{\tor}{\operatorname{Tor}}
\newcommand{\holim}{\operatornamewithlimits{holim}}
\newcommand{\Rlim}{\operatornamewithlimits{lim^1}}
\newcommand{\trc}{{\rm{trc}}}
\newcommand{\tr}{{\rm{tr}}}
\renewcommand{\leq}{\leqslant}
\renewcommand{\geq}{\geqslant}
\begin{document}
\title[On the algebraic K-theory of ku]{On the algebraic
K-theory\\ of the complex K-theory spectrum}
\author[Christian Ausoni]{Christian Ausoni}
\address{Mathematical Institute, University of Bonn, Germany}
\email{ausoni@math.uni-bonn.de}
\thanks{Published version DOI\,:~10.1007/s00222-010-0239-x.\\ 
Research supported in part by the Institute
Mittag-Leffler, Djursholm, and the Max Planck Institute for
Mathematics, Bonn.}
\begin{abstract}
Let $p\geq5$ be a prime, let $ku$ be the
connective complex $K$-theory spectrum,
and let $K(ku)$ be the algebraic $K$-theory spectrum of $ku$.
In this paper we study the $p$-primary homotopy
type of the spectrum $K(ku)$ by 
computing its mod $(p,v_1)$ homotopy groups. 
We show that up to a finite summand,
these groups form a finitely generated free module over the
polynomial algebra $\F_{p}[b]$,
where $b$ is a class of degree $2p+2$ defined as a ``higher Bott
element''.
\end{abstract}
\maketitle

\section{Introduction}
\label{kku-sec:intro}
The algebraic $K$-theory of a local or global number field $F$,
with suitable finite coefficients, is known 
to satisfy a form of Bott periodicity.
Bott periodicity refers here to the periodicity of
topological complex $K$-theory, and is an example of
$v_1$-periodicity in the sense of stable homotopy theory.
For example, if $p$ is an odd prime and if $F$ contains a primitive
$p$-th root of unity, then the mod $(p)$ algebraic $K$-theory
$K_*(F;\Z/p)$ of $F$ contains a non-nilpotent Bott element 
$\beta$ of degree $2$, with
\begin{equation*}
\beta^{p-1}=v_1\,.
\end{equation*}
In one of its reformulations~\cite{Rkku-DF85},\cite{Rkku-th85}, 
the Lichtenbaum-Quillen Conjecture
asserts that the localization
\begin{equation*}
K_*(F;\Z/p)\to K_*(F;\Z/p)[\beta^{-1}]
\end{equation*}
away from $\beta$ is an isomorphism in positive 
degrees. 
In particular, $K_*(F;\Z/p)$ is periodic of period~$2$ in
positive degrees. In the local case, this follows
from~\cite{Rkku-HM2}*{Theorem~D}. 
\par\smallskip
The $p$-local stable homotopy category also features higher forms of
periodicity~\cite{Rkku-HS98}, one for each integer $n\geq0$, referred to as
$v_n$-periodicity. It is detected for example by the $n$th Morava $K$-theory
$K(n)$, having coefficients $K(0)_*=\Q$ and
$K(n)_*=\F_p[v_n,v_n^{-1}]$ with
$|v_n|=2p^n-2$ if $n\geq1$. 
The study of $v_2$-periodicity is at the focus of current research in
algebraic topology, as illustrated for example by the efforts to
define the elliptic cohomology theory known as
topological modular forms~\cite{Rkku-Ho02}. 
\par
Waldhausen~\cite{Rkku-wal76} extended the definition of algebraic
$K$-theory to include specific ``rings up to homotopy'' 
called structured ring spectra, like $E_\infty$ ring spectra~\cite{Rkku-Ma77},
$S$-algebras~\cite{Rkku-EKMM}, or symmetric ring-spectra~\cite{Rkku-HSS}.
The chromatic red-shift conjecture~\cite{Rkku-ARg} of John
Rognes predicts that the algebraic $K$-theory
of a suitable $v_n$-periodic structured ring-spectrum is
essentially $v_{n+1}$-periodic, as illustrated above in the case of
number fields (which are $v_0$-periodic). 
For an example with the next level of periodicity, we
consider the algebraic $K$-theory of topological $K$-theory.
\par
\smallskip
Let $p\geq5$ be a prime, and let $ku_p$ denote the $p$-completed connective complex
$K$-theory spectrum with coefficients ${ku_p}_*=\Z_p[u]$,
$|u|=2$, where $\Z_p$ is the ring of $p$-adic integers.
Let $\ell_p$ be the Adams summand of $ku_p$ with
coefficients ${\ell_p}_*=\Z_p[v_1]$ and $v_1=u^{p-1}$.
In joint work with John Rognes~\cite{Rkku-AR02}, we have computed the
mod $(p, v_1)$ algebraic $K$-theory of the $S$-algebra
$\ell_p$, denoted
$V(1)_*K(\ell_p)$,
and we have shown that it is essentially
$v_2$-periodic. This computation provides a first example of
red-shift for non-ordinary rings.
\par
In this paper, following the discussion
in~\cite{Rkku-Au05}*{Section~10}, we interpret $ku_p$ as a tamely
ramified extension of $\ell_p$ of degree $p-1$, and we compute
$V(1)_*K(ku_p)$. As expected, the result is again
essentially periodic. However, 
$V(1)_*K(ku_p)$ has a shorter period\,: its periodicity is
given by multiplication with a higher Bott element
$b\in V(1)_*K(ku_p)$, of degree $2p+2$.
We defer a definition of $b$ to Section~\ref{kku-sec:units} below, and 
summarize our main result in the following statement.
\begin{theorem}\label{kku-thm:main}
Let $p\geq 5$ be a prime. The higher Bott element
$b\in V(1)_{2p+2}K(ku_p)$
is non-nilpotent and satisfies the relation 
\begin{equation*}
b^{p-1}=-v_2\,.
\end{equation*}
Let $P(b)$  denote the polynomial
$\F_p$-sub-algebra of $V(1)_*K(ku_p)$ generated by $b$. Then there is a short
exact sequence of graded $P(b)$-modules 
\begin{equation*}
0\to \Sigma^{2p-3}\F_p \to V(1)_*K(ku_p) \to F \to 0\,,
\end{equation*}
where $\Sigma^{2p-3}\F_p$ is the
sub-module  of $b$-torsion elements and
$F$ is a free $P(b)$-module on
$8+4(p-1)$ generators.
\end{theorem}
\par
A detailed
description of the free $P(b)$-module $F$ is given in 
Theorem~\ref{kku-thm:kku}.
The proof is based on evaluating the cyclotomic trace
map~\cite{Rkku-BHM}
\begin{equation*}
{\trc}:K(ku_p)\to TC(ku_p)
\end{equation*}
to topological cyclic homology.
We emphasize that the higher Bott element $b$ is not the
reduction of a class in the mod $(p)$ or integral homotopy
of $K(ku_p)$.
\par
The cyclic subgroup $\Delta\subset \Z_p^\times$ of order
$p-1$ acts on $ku_p$ by 
$p$-adic Adams operations.
The Adams summand is defined as the homotopy fixed-point
spectrum $\ell_p = ku_p^{h\Delta}$, and
$\Delta$ qualifies as the Galois group of the tamely
ramified extension $\ell_p\to ku_p$ of commutative
$S$-algebras given by the inclusion
of homotopy fixed-points. 
We proved in~\cite{Rkku-Au05}*{Theorem~10.2} that the induced map
$K(\ell_p)\to K(ku_p)$ factors through a weak
equivalence
\begin{equation*}
  K(\ell_p)\stackrel{\simeq}{\lr} K(ku_p)^{h\Delta}
\end{equation*}
after $p$-completion.  
The mod $(p,v_1)$ homotopy groups of $K(\ell_p)$ and
$K(ku_p)$ are related as follows.
\begin{proposition}\label{kku-prop:desc}
Let $i_*:V(1)_*K(\ell_p)\to V(1)_*K(ku_p)$ be the
homomorphism induced by the extension
of $S$-algebras $\ell_p\to ku_p$.
\begin{itemize} 
\item[\textup{(a)}] The homomorphism $i_*$
factors through an isomorphism 
\begin{equation*}
V(1)_*K(\ell_p)\cong\big(V(1)_*K(ku_p)\big)^\Delta\subset V(1)_*K(ku_p)
\end{equation*}
onto the classes fixed by the Galois group.
The higher Bott element $b$ is not fixed under the action of
$\Delta$, but $b^{p-1}=-v_2$ is, accounting for the
$v_2$-periodicity of $V(1)_*K(\ell_p)$.
\item[\textup{(b)}]
The homomorphism
\begin{equation*}
\mu:P(b) \otimes_{P(v_2)} V(1)_*K(\ell_p) \to V(1)_*K(ku_p) 
\end{equation*}
induced by $i_*$ and the $P(b)$-action
has finite kernel and cokernel, and is an isomorphism in degrees
larger than $2p^2-4$. By localizing away from $b$, we obtain
an isomorphism of $P(b, b^{-1})$-modules
\begin{equation*}
P(b,b^{-1}) \otimes_{P(v_2)} V(1)_*K(\ell_p)
\stackrel{\cong}{\lr} 
V(1)_*K(ku_p)[b^{-1}]\,. 
\end{equation*}
\end{itemize}
\end{proposition}
\par
In particular,
the $P(b)$-module $V(1)_*K(ku_p)$ is almost the module obtained 
from the $P(v_2)$-module $V(1)_*K(\ell_p)$ by the extension 
$P(v_2)\subset P(b)$ of scalars. 
The kernel of $\mu$ consists of $b$-multiples of the $v_2$-torsion
elements, and we have a non-trivial cokernel because 
some of the $P(v_2)$-module generators of $V(1)_*K(\ell_p)$ are
multiples of $b$ in $V(1)_*K(ku_p)$, see
Corollary~\ref{kku-cor:lvsku}.
\par
Notice that for the cyclotomic extension
$\Z_p\to\Z_p[\zeta_p]$ of complete discrete valuation
rings with Galois group $\Delta$ (where $\zeta_p$ is a
primitive $p$th root of unity), 
we have corresponding results in mod $(p)$ algebraic $K$-theory.
In effect, the natural homomorphism
$K_*(\Z_p;\Z/p)\to K_*(\Z_p[\zeta_p];\Z/p)$ factors through
an isomorphism onto the $\Delta$-fixed classes.
The Bott class $\beta\in K_2(\Z_p[\zeta_p];\Z/p)$ 
is not fixed under $\Delta$,
but $\beta^{p-1}=v_1$ is. This accounts for the fact that
$K_*(\Z_p[\zeta_p];\Z/p)$ has a shorter period than
$K_*(\Z_p;\Z/p)$. Moreover, the $P(\beta)$-module
$K_*(\Z_p[\zeta_p];\Z/p)$ is essentially obtained
from the $P(v_1)$-module $K_*(\Z_p;\Z/p)$ by the
extension $P(v_1)\subset P(\beta)$ of scalars.
These facts are extracted from computations by Hesselholt and
Madsen~\cite{Rkku-HM2}*{Theorem~D}.
We therefore interpret Proposition~\ref{kku-prop:desc} as follows~:
up to a chromatic shift of one in the sense
of stable homotopy theory, the
algebraic $K$-theory spectra of the tamely ramified extensions
\begin{equation*}
\xymatrix{
\Z_p[\zeta_p]&\ar@{}[d]|{\txt{and}}&ku_p\\
\Z_p\ar[u]^\Delta&&\ell_p\ar[u]_\Delta
}
\end{equation*}
have a comparable formal structure.
\par
This example of red-shift provides evidence that
structural results for the algebraic $K$-theory of ordinary
rings might well be generalized to provide more conceptual 
descriptions of the algebraic $K$-theory of $S$-algebras.
See Remarks~\ref{kku-rem:algclosed} and~\ref{kku-rem:HMlog} for
a discussion of the results we have in mind here.
\par\smallskip
We now turn to the algebraic $K$-theory $K(ku)$ of 
the (non $p$-completed) connective complex $K$-theory spectrum
$ku$, with coefficients $ku_*=\Z[u]$, $|u|=2$. 
The $p$-completion $ku\to ku_p$ induces a map 
\begin{equation*}
\kappa: K(ku)\to K(ku_p)\,, 
\end{equation*}
and the higher Bott element $b\in V(1)_{2p+2}K(ku_p)$ is in fact
defined as the image of a class with same name in
$V(1)_{2p+2}K(ku)$.
The difference between $K(ku)$ and $K(ku_p)$ can be measured
by means of the homotopy Cartesian square after $p$-completion
\begin{equation*}
\xymatrix{
K(ku)\ar[d]_\kappa\ar[r]^\pi &
K(\Z)\ar[d]^\kappa     \\
K(ku_p)\ar[r]^\pi       &
K(\Z_p)
}
\end{equation*}
of Dundas~\cite{Rkku-Du}*{page~224}.
Here $\pi$ denotes the map induced in $K$-theory by the
zeroth Postnikov sections $ku\to H\Z$ and $ku_p\to H\Z_p$,
where $HR$ is the Eilenberg-Mac\,Lane spectrum of the ring
$R$. 
The homotopy type of the $p$-completion of $K(\Z_p)$ has been computed 
by B\"okstedt, Hesselholt and
Madsen~\cite{Rkku-HM1},\cite{Rkku-BoM95}. The
Lichtenbaum-Quillen Conjecture for $K(\Z)$ (see for 
example~\cite{Rkku-mit}*{\S6}) implies that the homotopy
fiber of $K(\Z)\to K(\Z_p)$ has finite $V(1)$-homotopy
groups, which are concentrated in degrees smaller than
$2p-1$.
This implies the result below.
In fact there seems to be some consensus that work 
of Vladimir Voevodsky and Markus Rost should
imply the Lichtenbaum-Quillen Conjecture, but to our
knowledge this has not appeared in written form. We therefore
keep it as an assumption in the following results.
\begin{proposition}\label{kku-prop:integral}
Let $p\geq5$ be a prime, and assume that the Lichtenbaum-Quillen 
Conjecture for $K(\Z)$ holds at $p$. Then the homomorphism of $P(b)$-modules
\begin{equation*}
\kappa_*:V(1)_*K(ku)\to V(1)_*K(ku_p)
\end{equation*}
is an isomorphism in degrees larger than $2p-1$.
Localizing the $V(1)$-homotopy groups away from $b$, we obtain an isomorphism
\begin{equation*}
  V(1)_*K(ku)[b^{-1}]\cong V(1)_*K(ku_p)[b^{-1}]
\end{equation*}
of $P(b,b^{-1})$-algebras.
\end{proposition}
\par
This result is of interest beyond algebraic $K$-theory.
Baas, Dundas and Rognes have proposed 
a geometric definition of a cohomology theory
derived from a suitable notion of bundles of complex two-vector
spaces~\cite{Rkku-BDR}.
These are a two-categorical analogue of the ordinary complex
vector bundles which enter in the geometric definition of
topological $K$-theory.
They conjectured in~\cite{Rkku-BDR}*{5.1} that the spectrum
representing
this new theory is weakly homotopy
equivalent to $K(ku)$, and this was proved by these
authors and Birgit Richter in~\cite{Rkku-BDRR2}.
The next statement follows from Theorem~\ref{kku-thm:main} and
Proposition~\ref{kku-prop:integral}.
\begin{proposition}\label{kku-prop:tele}
If the Lichtenbaum-Quillen Conjecture for $K(\Z)$ holds, 
then at any prime $p\geq5$ the spectrum $K(ku)$ is of telescopic
complexity two in the sense of~\cite{Rkku-BDR}*{6.1}.
\end{proposition}
\par
This result was anticipated
in~\cite{Rkku-BDR}*{\S6}, and ensures that the 
cohomology theory derived from two-vector bundles is, 
from the view-point of stable homotopy theory, a
legitimate candidate for elliptic cohomology.
\par\smallskip
The computations presented in this paper fail at the primes
$2$ and $3$, because of the non-existence of the ring-spectrum
$V(1)$. Theoretically, computations in mod $(p)$ homotopy or
in integral homotopy could also be carried out, but the 
algebra seems quite intractable. Another 
approach~\cite{Rkku-BR05},\cite{Rkku-LN05} is via
homology computations. There are ongoing projects in this
direction by Robert Bruner, Sverre Lun{\o}e-Nielsen and John Rognes. 
\par\smallskip
Up to degree three, the integral homotopy groups of $K(ku)$ can be
computed essentially by using the map $\pi:K(ku)\to K(\Z)$
introduced above.  The map
$\pi_*:K_*(ku)\to K_*(\Z)$ is
$3$-connected, so that 
\begin{equation*}
K_0(ku)\cong \Z,\ K_1(ku)\cong \Z/2\ \textup{ and }\
K_2(ku)\cong\Z/2\,.
\end{equation*}
Here $K_1(ku)$ and $K_2(ku)$ are generated by
the image of $\eta\in \pi_1S$ and $\eta^2\in \pi_2S$,
respectively, under
the unit $S\to K(ku)$.
Let $w:BBU_\otimes \to \Omega^\infty K(ku)$ be the map induced by the
inclusion of units, see~\eqref{kku-eq:unit1}.
There is a non-split extension
\begin{equation*}
0\to
\pi_3(BBU_\otimes)\stackrel{w_*}{\lr}K_3(ku)\stackrel{\pi_*}{\lr}K_3(\Z)\to 0
\end{equation*}
with $\pi_3(BBU_\otimes)\cong\Z\{\mu\}$, $K_3(ku)\cong
\Z\{\varsigma\}\oplus\Z/24\{\nu\}$ and
$K_3(\Z)\cong\Z/48\{\lambda\}$, where $\nu$ is the image of
the Hopf class $\nu$, which generates $\pi_3 S\cong\Z/24$. 
We have $w_*(\mu)=2\varsigma-\nu$ and
$\pi_*(\varsigma)=\lambda$. See~\cite{Rkku-ADR08} for details.
This indicates that the integral homotopy groups $K_*(ku)$ contain 
intriguing non-trivial extensions from subgroups in
$\pi_*S$, $\pi_* BBU_\otimes$ and $K_*(\Z)$. 
\par\smallskip
The rational algebraic $K$-groups of $ku$ are well understood.
In joint work with John Rognes~\cite{Rkku-ARQ}, we have proved that
after rationalization, the sequence
\begin{equation*}
BBU_\otimes\stackrel{w}{\lr}\Omega^\infty K(ku)
\stackrel{\pi}{\lr} \Omega^\infty K(\Z)
\end{equation*}
is a split homotopy fibre-sequence. A rational splitting of
$w$ is provided by a rational determinant map $\Omega^\infty K(ku)\to
(BBU_\otimes)_\Q$. 
In particular, by Borel's computation~\cite{Rkku-Bo74} of
$K_*(\Z)\otimes\Q$, there is a rational equivalence
\begin{equation*}
\Omega^\infty K(ku)\simeq_\Q SU\times (SU/SO)\times\Z\,.
\end{equation*}
All but finitely many of the non-torsion classes in the integral 
homotopy groups $\pi_*K(ku)$ detected by this equivalence 
reduce mod $(p)$ to multiples of $v_1$, and hence
reduce to zero in $V(1)_*K(ku)$. 
\par\smallskip
We briefly discuss the contents of this paper. In
Section~\ref{kku-sec:kz3}, 
we study the $V(1)$-homotopy of the Eilenberg-Mac\,Lane space 
$K(\Z,3)$, which is a subspace of the space of units of $ku$. 
In Section~\ref{kku-sec:units}, we define low-dimensional classes in 
$V(1)_*K(ku)$ corresponding to units of $ku$, and in
particular we
introduce the higher Bott element.
We prove in Section~\ref{kku-sec:trace} that these classes are non-zero by
means of the B\"okstedt trace map
\begin{equation*}
\tr:K(ku)\to \thh(ku)
\end{equation*} 
to topological Hochschild homology.
In Section~\ref{kku-sec:low}, we compute $V(1)_nK(ku_p)$ for $n\leq
2p-2$. This complements the computations in higher degrees 
provided by the cyclotomic trace
\begin{equation*}
\trc:K(ku_p)\to \tc(ku_p)
\end{equation*}
to topological cyclic homology.
In Section~\ref{kku-sec:fixed} we compute the various homotopy fixed points
of $\thh(ku_p)$ under the action of the cyclic groups $C_{p^n}$ and the
circle, which are the ingredients for the computation
of $V(1)_*TC(ku_p)$ in Section~\ref{kku-sec:tc} . 
In Section~\ref{kku-sec:algK} we prove Theorem~\ref{kku-thm:main} 
on the structure of $V(1)_*K(ku_p)$ stated above. 
We also give a 
computation of $V(1)_*K(KU_p)$ for $KU_p$ the $p$-completed periodic
$K$-theory spectrum, up to some indeterminacy.
\begin{notations}
Throughout the paper, unless stated otherwise, $p$ will be a 
fixed prime with $p\geq5$, and
$\Z_p$ will denote the $p$-adic integers.
For an $\F_p$-vector space $V$, let $E(V)$, $P(V)$ 
and $\Gamma(V)$ be the exterior algebra,
polynomial algebra and divided power algebra on $V$, 
respectively. If $V$ has a basis 
$\{x_1,\dots,x_n\}$, we write $V=\F_p\{x_1,\dots,x_n\}$
and $E(x_1,\dots,x_n)$, 
$P(x_1,\dots,x_n)$ and $\Gamma(x_1,\dots,x_n)$
for these algebras. 
By definition, $\Gamma(x)$ is the $\F_p$-vector space 
$\F_p\{\gamma_kx\,|\,k\geq0\}$ with 
product given by $\gamma_ix\cdot\gamma_jx=\binom{i+j}i\gamma_{i+j}x$, 
where $\gamma_0x=1$ and $\gamma_1x=x$.
Let $P_h(x)=P(x)/(x^h)$ be the 
truncated polynomial algebra of height $h$.
For an algebra $A$, we denote by $A\{x_1,\dots,x_n\}$ 
the free $A$-module generated by 
$x_1,\dots,x_n$.
\par
If $Y$ is a space and $E_*$ is a homology theory, 
such as mod $(p)$ homology, $V(1)$-homotopy or Morava
$K$-theory $K(2)_*$, we denote by $E_*(Y)$ the 
unreduced $E_*$-ho\-mo\-lo\-gy of $Y$, which we identify with
the $E_*$-homology of the suspension
spectrum $\Sigma^\infty(Y_+)$, where $Y_+$ denotes $Y$ with
a disjoint base-point added. We usually write
$\Sigma^\infty_+Y$ instead of $\Sigma^\infty(Y_+)$. 
\par
The reduced $E_*$-homology of a pointed space $X$ is
denoted $\wti E_*(X)$. 
We denote $\pi_*X$ the (unstable) homotopy groups of
$X$, and $\pi_*\Sigma^\infty X$ its stable
homotopy groups. 
\par
If $f:A\to B$ is a map of $S$-algebras, we also denote by~$f$ 
its image under various functors like $\thh$,
$\tc$ or $K$.
\par
In our computations with spectral sequences, we often determine
a differential $d$ only up to multi\-pli\-ca\-tion by a unit.
We use the notation $d(x)\doteq y$ to indicate that the equation
$d(x)=\alpha y$ holds for some unit $\alpha\in\F_p$.
Classes surviving to the $E^r$-term of a spectral sequence, for $r\geq 3$,
are often given as a product of classes in the $E^2$-term.
To improve the readability, we denote the product of two
classes $x,y$ in $E^r$ by $x\cdot y$. 
\end{notations}
\section {On the $V(1)$-homotopy of $K(\Z,3)$}
\label{kku-sec:kz3}
If $G$ is a
topological monoid, let us denote by $BG$ its classifying
space, obtained by realization of the bar construction, see
for example~\cite{Rkku-RW}*{\S1}. If $G$ is an Abelian
topological group, then so is $BG$.  The space $BG$
is equipped with the bar filtration
\begin{equation}\label{kku-eq:barfiltr}
\{*\}=B_0\subset B_1\subset B_2\subset\dots\subset
B_{n-1}\subset B_{n}\subset\dots BG\,, 
\end{equation}
with filtration quotients $B_{n}/B_{n-1}\cong \Sigma^n
(G^{\wedge n})$. 
In particular, we have a map 
\begin{equation}\label{kku-eq:1squelet}
s:\Sigma G=B_1\subset BG\,,
\end{equation}
which in any homology theory $E_*$ induces a map
\begin{equation*}
\sigma:E_{*}G\to E_{*+1}BG
\end{equation*}
called the {\sl suspension}.
If $E_*$ is a multiplicative homology theory satisfying the
K\"un\-neth isomorphism, we have the bar spectral
sequence~\cite{Rkku-RW}*{\S2} 
\begin{equation*}
\begin{split}
  E^1_{s,*}(G)  &  =\wti E_*(G)^{\otimes_{E_*}s},\\
E^2_{s,t}(G)  &  =\tor^{E_*(G)}_{s,t}(E_*,E_*)
\Rightarrow E_{s+t}(BG)
\end{split}
\end{equation*}
associated to the bar filtration~\eqref{kku-eq:barfiltr}.
\par
Let $K(\Z,0)$ be equal to $\Z$ as a discrete topological group, and
for $m\geq1$, we define recursively the Eilenberg-Mac\,Lane space
$K(\Z,m)$ as the Abelian topological group $BK(\Z,m-1)$. 
We recall Cartan's computation of the algebra $H_*(K(\Z,m);\F_p)$ 
for $p$ an odd prime and $m=2,3$.
The generators are
constructed explicitly from the unit $1\in
H_*(K(\Z,0);\F_p)$ by means of the suspension $\sigma$ and
two further operators 
\begin{equation*}
\begin{split} 
\varphi: H_{2q}(K(\Z,m);\F_p) & \to
H_{2pq+2}(K(\Z,m+1);\F_p)\ \ {\textup{and}}\\
\gamma_p: H_{2q}(K(\Z,m);\F_p) & \to
H_{2pq}(K(\Z,m);\F_p)\,,
\end{split}
\end{equation*}
called  the {\sl
transpotence}~\cite{Rkku-Car}*{page~6-06} and the
{\sl $p$-th divided power}~\cite{Rkku-Car}*{page~7-07}, respectively.
The transpotence is an additive
homomorphism since $p$ is odd. 
For $x\in H_{2q}(K(\Z,m);\F_p)$, the class $\varphi(x)$ is
represented, for example, by
\begin{equation*}
x^{p-1}\otimes x\in E^1_{2,2pq}\big(K(\Z,m)\big)
\end{equation*}
in the bar spectral sequence. The algebra
$H_{*}(K(\Z,m);\F_p)$ has the structure of an algebra with
divided powers, which are uniquely determined by $\gamma_p$.
\begin{theorem}[Cartan]\label{kku-thm:cartan}
Let $p$ be an odd prime. 
There are isomorphisms of $\F_p$-alge\-bras with divided
powers
\begin{equation*}
\Gamma(y)\stackrel{\cong}{\lr}
H_*(K(\Z,2);\F_p)\,
\end{equation*}
given by $y \mapsto \sigma\sigma(1)$, with $|y|=2$, and 
\begin{equation*}
\bigotimes_{k\geq0}E(e_k)\otimes\Gamma(f_k)\stackrel{\cong}{\lr}
H_*(K(\Z,3);\F_p)\,,
\end{equation*}
given by 
$e_k \mapsto \sigma\gamma_p^k\sigma\sigma(1)$ and
$f_k \mapsto \varphi\gamma_p^k\sigma\sigma(1)$, with
degrees $|e_k|=2p^k+1$ and $|f_k|= 2p^{k+1}+2$. 
For $k\geq0$, the generators $f_k$ and
$e_{k+1}$ are related by a primary mod $(p)$ homology Bockstein
\begin{equation*}
\beta(f_k)=e_{k+1}\,.
\end{equation*}
\end{theorem}
\begin{proof}
The computation of $H_*(K(\Z,m);\F_p)$ as an algebra is
given in~\cite{Rkku-Car}*{Th\'eo\-r\`eme fondamental, p. 9-03}.
The Bockstein
relation $\beta(f_k)=e_{k+1}$ is established
in~\cite{Rkku-Car}*{page~8-04}.
\end{proof}
\par
Ravenel and Wilson~\cite{Rkku-RW} make use of the bar spectral sequence to compute the
Morava $K$-theory $K(n)_*K(\pi, m)$ as an algebra when $\pi =\Z$ or $\Z/p^j$. 
All generators 
can be defined explicitly, starting with
the unit $1\in K(n)_*K(\pi,0)$ and using the suspension,
divided powers, transpotence and the
Hopf-ring structure on $K(n)_*K(\pi,*)$.  
We refer to~\cite{Rkku-RW}*{5.6 and 12.1} for the
following result, and for the definition of the generators
$\beta_{(k)}$ and $b_{(2k,1)}$. 
\begin{theorem}[Ravenel-Wilson]\label{kku-thm:RW}
Let $p\geq3$ be a prime and let $K(2)$ be the Morava
K-theory spectrum with coefficients $K(2)_*=\F_p[v_2,v_2^{-1}]$.
There are isomorphisms of $K(2)_*$-alge\-bras
\begin{equation*}
K(2)_*K(\Z,2)\cong
K(2)_*[\beta_{(k)}\,|\,k\geq0]/(\beta_{(0)}^p,
\beta_{(k+1)}^p-v_2^{p^k}\beta_{(k)}\,|\,k\geq
0)
\end{equation*}
where $|\beta_{(k)}|=2p^k$, and 
\begin{equation*}
K(2)_*K(\Z,3)\cong
K(2)_*[b_{(2k,1)}\,|\,k\geq0]/(b_{(2k,1)}^p+v_2^{p^k}b_{(2k,1)}\,|\,k\geq0)
\end{equation*}
where $|b_{(2k,1)}|=2p^k(p+1)$. 
The class $\beta_{(0)}\in K(2)_2K(\Z,2)$ is equal to $\sigma\sigma(1)$, and
the class $b_{(0,1)}\in K(2)_{2p+2}K(\Z,3)$ is the
transpotence of $\beta_{(0)}$. 
\end{theorem}
\par
We now turn to $V(1)$-homotopy. 
For an integer $n\geq0$, we denote by $V(n)$ the Smith-Toda complex~\cite{Rkku-To71}, 
with mod $(p)$ homology given by 
\begin{equation*}
H_*(V(n);\F_p)\cong E(\tau_0,\dots,\tau_n)
\end{equation*} 
as a left sub-comodule of the dual Steenrod algebra. 
In particular, $V(0)=S/p$ is the mod $(p)$ Moore spectrum,
and the spectra
$V(0)$ and $V(1)$ fit in cofibre sequences
\begin{equation*}
S\stackrel{p}{\lr} S\stackrel{i_0}{\lr} V(0)
\stackrel{j_0}{\lr} \Sigma S
\end{equation*}
and
\begin{equation*}
\Sigma^{2p-2}V(0)\stackrel{v_1}{\lr} V(0)\stackrel{i_1}{\lr} V(1)
\stackrel{j_1}{\lr}
\Sigma^{2p-1}V(0)\,,
\end{equation*}
where $v_1$ is a periodic map. 
For $n=0,1$ and $p\geq 5$, the spectrum $V(n)$ is a commutative 
ring spectrum~\cite{Rkku-Ok}, and its ring of coefficients $V(n)_*$ 
is an $\F_p$-algebra which contains a non-nilpotent class
$v_{n+1}$, of degree $2p^{n+1}-2$. 
We call ``$V(n)$-homotopy'' the homology theory associated to
the spectrum $V(n)$. In other words, the 
$V(n)$-homotopy groups of a spectrum $X$ are defined by
\begin{equation*}
V(n)_*X=\pi_*(V(n)\wedge X)\,.
\end{equation*}
Notice that $V(0)_*X$ is denoted $\pi_*(X;\Z/p)$ by some
authors, and called the mod $(p)$ homotopy groups of $X$.  
By analogy, we sometimes call $V(1)_*X$ the mod $(p,v_1)$
homotopy groups of $X$.
If $Y$ is a space, then $V(n)_*Y$ is defined as
$V(n)_*\Sigma^\infty_+ Y$.
\par
The primary mod $(p)$ homotopy
Bockstein $\beta_{0,1}:V(0)_*X\to V(0)_{*-1}X$ is the
homomorphism induced by $(\Sigma i_0)j_0$, and 
the primary mod $(v_1)$ homotopy Bockstein
$\beta_{1,1}:V(1)_*X\to V(1)_{*-2p+1}X$ is the homomorphism
induced by $(\Sigma^{2p-1}i_1)j_1$.
The homomorphisms ${i_0}_*:\pi_*(X)\to V(0)_*X$ and
${i_1}_*:V(0)_*X\to V(1)_*X$ are called the mod $(p)$
reduction and the mod $(v_1)$ reduction, respectively.
\par
Let $H\F_p$ be the Eilenberg-Mac\,Lane spectrum of $\F_p$.
The unit map $S\to H\F_p$ factors through a map of ring
spectra $h:V(1)\to H\F_p$, which induces an injective
homomorphism in mod $(p)$ homology. Identifying the homology
of $V(1)$ with its image in the dual Steenrod algebra $A_*$,
we obtain the isomorphism 
\begin{equation*}
H_*(V(1);\F_p)\cong E(\tau_0,\tau_1)
\end{equation*}
of left $A_*$-comodule algebras mentioned above. 
Toda~\cite{Rkku-To71}*{Theorem~5.2} computed $V(1)_*$ in a range of
degrees for which the Adams spectral sequence collapses.
Up to some renaming of the classes, we deduce from his
theorem that for $p\geq5$ there is an isomorphism of
$P(v_2)\otimes P(\beta_1)$-modules
\begin{equation}\label{kku-eq:toda}
P(v_2)\otimes
P(\beta_1)\otimes\F_p\{1,\alpha_1,\beta_1',(\alpha_1\beta_1)^\sharp \}
\to V(1)_*
\end{equation}
in degrees $*< 4p^2-2p-4$. 
The classes $\alpha_1$ and $\beta_1$ are the mod $(p,v_1)$
reduction of the classes with same name in $\pi_*(S)$, of
degrees $2p-3$ and $2p^2-2p-2$, respectively. The class
$\beta_1'$ is the mod $(v_1)$ reduction of the class with same name in
$V(0)_*$ that supports a primary mod $(p)$ homotopy Bockstein
$\beta_{0,1}(\beta_1')=\beta_1$, and is of degree
$2p^2-2p-1$. 
The classes $v_2$ and $(\alpha_1\beta_1)^{\sharp}$, of degree $2p^2-2$
and $2p^2+2p-6$ respectively, support a primary mod $(v_1)$ homotopy
Bockstein, given by $\beta_{1,1}(v_2)=\beta_1'$ and
$\beta_{1,1}((\alpha_1\beta_1)^\sharp)=\alpha_1\beta_1$.
The class $v_2$ is non-nilpotent.
The lowest-degree class in $V(1)_*$ that is not in the image
of~\eqref{kku-eq:toda} is the mod
$(p,v_1)$ reduction of the class $\beta_2$ in $\pi_*(S)$, of
degree $4p^2-2p-4$.
\par
If $X$ is a connective spectrum of finite type, the
Atiyah-Hirzebruch spectral sequence
\begin{equation}\label{kku-eq:AHss}
E^2_{s,t}=H_s(X;\F_p)\otimes V(1)_t \Rightarrow V(1)_{s+t}X
\end{equation} 
converges strongly, and we can use it to compute $V(1)_*X$ in low degrees. 
The first non-trivial Postnikov invariant of $V(1)$ 
is Steenrod's reduced power operation $P^1$,
corresponding to the first possibly non-trivial differential
of the spectral sequence on the zeroth line, see Remark~\ref{kku-rem:AHss}. 
This operation detects the class $\alpha_1$, which belongs to the
kernel of the Hurewicz homomorphism $V(1)_*\to
H_*(V(1);\F_p)$. In some more details,
we have a commutative diagram
\begin{equation}\label{kku-eq:postnikov}
\xymatrix{
& 
V(1)\ar[d]^\rho \ar@/^/[dr]^h  &
&  \\
\Sigma^{2p-3}H\F_p \ar[r]^-g & 
V(1)[2p-3] \ar[r]^-h         & 
H\F_p \ar[r]^-{P^1}      &
\Sigma^{2p-2}H\F_p\,,
}
\end{equation}
where $\rho$ is the $(2p-3)$th-Postnikov section, and
the horizontal sequence is a cofibre sequence.
Notice that by~\eqref{kku-eq:toda} the map $\rho$ is $(2p^2-2p-2)$-connected, so
that under our assumptions on $X$ we have a well defined homomorphism 
\begin{equation*}
\alpha=(\rho_*)^{-1}g_* :H_{n-2p+3}(X;\F_p)\to V(1)_nX
\end{equation*}
for $n\leq 2p^2-2p-3$.
\begin{lemma}\label{kku-lem:approx}
Let $X$ be a connective spectrum of finite type,
and let $p\geq3$ be a prime.
For $n\leq2p^2-2p-3$, the group $V(1)_nX$ fits in
an exact sequence
\begin{equation*}
\begin{split}
H_{n+1}(X;\F_p) & 
\stackrel{(P^1)^*}{\lr}
H_{n-2p+3}(X;\F_p)\stackrel{\alpha}{\lr}
V(1)_nX \stackrel{h_*}{\lr}\\
H_n(X;\F_p)&\stackrel{(P^1)^*}\lr H_{n-2p+2}(X;\F_p)\,.
\end{split}
\end{equation*}
Here $(P^1)^*$ denotes the homology operation dual to $P^1$.
If $X$ is a ring spectrum then $\alpha$ sends the unit $1\in
H_0(X;\F_p)$ to $\alpha_1$. 
Moreover, for any $X$ and any $n\geq0$, we have a commutative diagram
\begin{equation*}
\xymatrix{
V(1)_nX \ar[r]^-{h_*} \ar[d]_{\beta_{1,1}} &
H_n(X;\F_p) \ar[d]^{Q_1^*}      \\
V(1)_{n-2p+1}X \ar[r]^-{h_*}               &
H_{n-2p+1}(X;\F_p)
}
\end{equation*}
relating the primary mod $(v_1)$ homotopy Bockstein
$\beta_{1,1}$ to the homology operation 
$Q_1^*$ dual to Milnor's primitive $Q_1=P^1\delta-\delta P^1\in A$.
\end{lemma}
\begin{proof}
This exact sequence is the sequence associated to the
cofibre sequence in \eqref{kku-eq:postnikov}, where we have replaced
$V(1)[2p-3]_nX$ by $V(1)_nX$ via $\rho_*$, which is an
isomorphism for these values of $n$, by strong convergence
of the Atiyah-Hirzebruch spectral sequence. 
The assertion on $\alpha_1$ is true if $X=S$, and 
follows by naturality for $X$ an arbitrary ring spectrum.
\par
The self-map
$f=(\Sigma^{2p-1}i_1)j_1$ of $V(1)$, which induces
$\beta_{1,1}$, is
given in mod $(p)$ homology by the homomorphism $f_*:E(\tau_0,\tau_1)\to
E(\tau_0,\tau_1)$ of degree $1-2p$ with  $f_*(1)=f_*(\tau_0)=0$,  
$f_*(\tau_1)=1$ and
$f_*(\tau_0\tau_1)=\tau_0$. We have a commutative
diagram
\begin{equation*}
\xymatrix{
V(1)_*X\ar[r]^-{\beta_{1,1}} \ar[d]_{g_*} &
V(1)_*X \ar[d]_{g_*}     \\
E(\tau_0,\tau_1)\otimes
H_*(X;\F_p)\ar[r]^-{f_*\otimes1}\ar[d]_{e_*} &
E(\tau_0,\tau_1)\otimes H_*(X;\F_p)\ar[d]^-{\mu} \\ 
A_*\otimes H_*(X;\F_p)\ar[r]^-{{\tau_1}^*\otimes1}& 
H_*(X;\F_p)
\,.
}
\end{equation*}
The horizontal arrows are of degree $1-2p$, and
${\tau_1}^*:A_*\to \F_p$ is the dual of
$\tau_1$ with respect to the
standard basis $\{\tau(E)\xi(R)\}$ of $A_*$ given in~\cite{Rkku-Mil}*{\S 6}. 
The homomorphism $g_*$ is induced in homotopy by the smash product of the
unit $S\to H\F_p$ with the
identity of $V(1)\wedge X$,  $\mu$ is
induced by the right homotopy action $H\F_p\wedge V(1)\to
H\F_p$, and $e_*$ is induced by $1\wedge
h\wedge1:H\F_p\wedge V(1)\wedge X\to H\F_p\wedge H\F_p\wedge
X$. 
We have $\mu g_*=h_*$ and $e_* g_*=\nu_* h_*$, where $\nu_*$
is the left $A_*$-coaction on $H_*(X;\F_p)$. 
This completes the proof since
$({\tau_1}^*\otimes1)\nu_*=Q_1^*$ by definition of $Q_1$,
see~\cite{Rkku-Mil}*{page~163}.
\end{proof}
\begin{remark}\label{kku-rem:AHss}
For $X$ connective, the Atiyah-Hirzebruch spectral
sequence~\eqref{kku-eq:AHss}
has only two non-trivial lines in internal degrees $t\leq
2p^2-2p-3$, corresponding to $1$ and $\alpha_1$ in $V(1)_*$,
see~\textup{\eqref{kku-eq:toda}}. 
The argument above shows that there is a differential 
\begin{equation*}
d^{2p-2}(z)=(P^1)^*(z)\alpha_1
\end{equation*}
for $z\in E^2_{*,0}$. 
In total degrees less than $2p^2-2p-3$ this is the only possibly
non-trivial differential.
\end{remark}
\begin{lemma}\label{kku-lem:v1kz2}
The map
\begin{equation*}
\F_p\{\alpha_1\}\oplus{P_p}(x)\to {V(1)}_*K(\Z,2)
\end{equation*}
given by $x\mapsto \sigma\sigma(1)$ with $|x|=2$ is an isomorphism in
degrees less than $4p-3$.
\end{lemma}
\begin{proof}
This follows from Theorem~\ref{kku-thm:cartan}, 
Lemma~\ref{kku-lem:approx} and the relation
\begin{equation}\label{kku-eq:p1}
(P^1)^*\big(\gamma_{k+p-1}(y)\big)=k\gamma_{k}(y)
\end{equation}
in $H_*(K(\Z,2);\F_p)\cong\Gamma(y)$. 
\end{proof}
\par
Consider the cofibration 
\begin{equation*}
B_1=\Sigma K(\Z,2)\stackrel{i}{\lr} B_2\stackrel{j}{\lr}
\Sigma^2(K(\Z,2)^{\wedge2})\to \Sigma^2  K(\Z,2)
\end{equation*}
extracted from the bar filtration~\eqref{kku-eq:barfiltr} of $K(\Z,3)$. It induces
an exact sequence
\begin{equation*}
V(1)_{*}\Sigma K(\Z,2)\stackrel{i_*}{\lr}
V(1)_*B_2\stackrel{j_*}{\lr}
\wti{V(1)}_{*}\Sigma^2(K(\Z,2)^{\wedge2})\stackrel{\Sigma^2\mu_*}{\lr}
V(1)_{*}\Sigma^2 K(\Z,2)\,,
\end{equation*}
where $\mu_*$ is induced by the product on $K(\Z,2)$.
We know  that $V(1)_{2p+1}K(\Z,2)=0$, by
Lemma~\ref{kku-lem:v1kz2},
which implies that the homomorphism
\begin{equation*}
V(1)_{2p+2}B_2\stackrel{j_*}{\lr}\wti{V(1)}_{2p}K(\Z,2)^{\wedge2}
\end{equation*}
is injective.
We know as well that the composition
\begin{equation*}
\wti{V(1)}_*K(\Z,2)\otimes
\wti{V(1)}_*K(\Z,2) \stackrel{k}{\lr}
\wti{V(1)}_{*}K(\Z,2)^{\wedge2} \stackrel{\mu_*}{\lr}
V(1)_*K(\Z,2)
\end{equation*}
sends the class $x^{p-1}\otimes x$ to zero. 
In particular, the class $k(x^{p-1}\otimes x)$ is in the image of $j_*$.
Let $\ti b'\in
V(1)_{2p+2}B_2$ be the unique class which satisfies the
equation
\begin{equation*}
j_*(\ti b')=\Sigma^2 k(x^{p-1}\otimes x)\,. 
\end{equation*}
\begin{definition}\label{kku-def:prebott}
We define the fundamental class $e'_0\in V(1)_3K(\Z,3)$ as the image
of the unit $1\in V(1)_0K(\Z,0)$ under the iterated
suspension $\sigma^3$.
We define 
\begin{equation*}
b'\in V(1)_{2p+2}K(\Z,3)
\end{equation*} 
as $b'={l_2}_*(\ti b')$,
where $l_2:B_2\to K(\Z,3)$ is the inclusion of the second
subspace in the bar filtration.
\end{definition}
\par
Notice that the definition of $b'$ in $V(1)$-homotopy, using $x^{p-1}\otimes
x$ as above, lifts the definition of the transpotence in the
homology of the bar construction. We use this fact in the
proof of the following proposition.
\begin{proposition}\label{kku-prop:prebott}
The class $b'\in V(1)_*K(\Z,3)$ is
non-nilpotent, and satisfies the relation
\begin{equation*}
b'^{p} = -v_2b'\,.
\end{equation*}
There is a primary mod $(v_1)$ homotopy Bockstein 
\begin{equation*}
\beta_{1,1}(b') = e'_0\,.
\end{equation*}
\end{proposition}
\begin{proof}
First, we notice that the $\F_p$-vector space
$V(1)_{2p^2+2p}K(\Z,3)$, which contains $b'^{p}$, is of rank
at most one. Indeed, 
consider the Atiyah-Hirzebruch spectral sequence
\begin{equation*}
E^2_{s,t}\cong H_s(K(\Z,3);\F_p)\otimes V(1)_t\Rightarrow
V(1)_{s+t}K(\Z,3)\,.
\end{equation*}
From Theorem~\ref{kku-thm:cartan} and the formula~\eqref{kku-eq:toda} 
for $V(1)_*$ in low degrees we deduce that
$E^2_{*,*}$ consists of
$\F_p\{f_0\cdot v_2,\ e_0\cdot f_0\cdot \alpha_1\cdot \beta_1\}$ 
in total degree $2p^2+2p$. 
Suspending the relation~\eqref{kku-eq:p1} for $k=1$ we get a
relation
\begin{equation}\label{kku-eq:p1e1}
(P^1)^*(e_1)=e_0\,.
\end{equation}
Notice that for degree reasons the class
$e_1\cdot f_0\cdot\beta_1\in E^2_{*,*}$ survives to
$E^{2p-2}_{*,*}$ as a product of $e_1$ and $f_0\cdot \beta_1$.
By Remark~\ref{kku-rem:AHss}, and since $f_0\cdot \beta_1$ is a cycle,
we have a differential
\begin{equation*} 
d^{2p-2}(e_1\cdot f_0\cdot \beta_1)=
e_0\cdot f_0\cdot \alpha_1\cdot \beta_1\,,
\end{equation*}
and this implies the claim on $V(1)_{2p^2+2p}K(\Z,3)$.
\par
The unit map $S\to K(2)$ factors through a map of ring
spectra $V(1)\to K(2)$. The induced ring homomorphism
\begin{equation*}
V(1)_*K(\Z,2)\to K(2)_*K(\Z,2) 
\end{equation*}
maps $x$ to $\beta_{(0)}$, since these classes are defined
as the double suspension of the unit in $V(1)_0K(\Z,0)$,
respectively $K(2)_0K(\Z,0)$. By construction, the
class $b'$ maps to the transpotence of $\beta_{(0)}$, which is
$b_{(0,1)}$. 
We deduce that the sub-$V(1)_*$-algebra of $V(1)_*K(\Z,3)$
generated by $b'$ maps surjectively onto the subalgebra
\begin{equation*}
P(v_2,b_{(0,1)})/(b_{(0,1)}^p+v_2b_{(0,1)})
\end{equation*}
of $K(2)_*K(\Z,3)$ generated by $v_2$ and $b_{(0,1)}$.
In particular $b'$ is non-nilpotent. Thus
$V(1)_{2p^2+2p}K(\Z,3)$ is of rank one, and injects
into $K(2)_{2p^2+2p}K(\Z,3)$. This implies the identity
$b'^{p} = -v_2b'$.
\par
To prove the Bockstein relation, we map to homology. The
Hurewicz homomorphism $h_*:V(1)_*K(\Z,3)\to H_*(K(\Z,3);\F_p)$ is an
isomorphism in degrees $3$ and $2p+2$, mapping
$e'_0$ to $e_0$ and $b'$ to the transpotence
$\varphi(y)=f_0$ of $y$.
We have a primary homology Bockstein
$\beta(f_0)=e_1$ by Theorem~\ref{kku-thm:cartan}, 
and combining with~\eqref{kku-eq:p1e1} 
we obtain $ (P^1)^*\beta(f_0)=e_0$. 
We also have $\beta(P^1)^*(f_0)=0$ for degree reasons.  
Finally,
\begin{equation*}
Q_1^*(f_0)=\big((P^1)^*\beta-\beta
(P^1)^*\big)(f_0)=e_0\,,
\end{equation*}
so by Lemma~\ref{kku-lem:approx} the relation
$\beta_{1,1}(b')=e'_0$ holds.
\end{proof}
\section{The units of $ku$ and the higher Bott element}
\label{kku-sec:units}
The aim of this section is to define low-dimensional classes in 
$V(1)_*K(ku)$ by using the inclusion of units.
\par
We recall from~\cite{Rkku-Ma77} or \cite{Rkku-Ma09}*{Definition~7.6} 
that the space of units $GL_1(A)$ 
of an $E_\infty$-ring spectrum $A$ is defined by
the following pull-back square of spaces
\begin{equation*}
\xymatrix{
GL_1(A) \ar[r] \ar[d]          &
\Omega^\infty A \ar[d]^{\pi_0} \\
GL_1(\pi_0A) \ar[r]           &
\ \pi_0A\,.
}
\end{equation*}
Taking the vertical fiber over $1\in GL_1(\pi_0A)$, we obtain a
fiber sequence of group-like $E_\infty$-spaces or infinite loop spaces
\begin{equation*}
SL_1(A)\to GL_1(A)\to GL_1(\pi_0A)\,,
\end{equation*}
with products given by the multiplicative structure of
$A$. Here we can assume that we have a model of $GL_1(A)$ 
and of $SL_1(A)$ which is actually a topological monoid, see for
example~\cite{Rkku-Sch04}*{\S 2.3}.
The functor $GL_1$ from $E_\infty$-ring spectra to
infinite loop spaces is right adjoint, up to homotopy,
to the suspension functor $\Sigma_+^\infty$. This follows
from~\cite{Rkku-Ma09}*{Lemma~9.6}.
\par
In the case of $ku$, the space $SL_1(ku)$ is commonly
denoted $BU_\otimes$.
This notation refers to the product of the underlying
$H$-space of $BU_\otimes$, which represents the tensor 
product of virtual line bundles.
\par
The first Postnikov section $\pi:BU_\otimes\to K(\Z,2)$, with
homotopy fiber denoted by $BSU_\otimes$, admits a section
$j:K(\Z,2)\simeq BU(1)\to BU_\otimes$. 
Here the map $j$ represents viewing a
line bundle as a virtual line bundle.   
Both $\pi$ and $j$
are infinite loop maps, and we have a splitting of infinite
loop-spaces
\begin{equation*}
BU_\otimes\simeq K(\Z,2)\times BSU_\otimes\,,
\end{equation*}
see~\cite{Rkku-Ma77}*{V.3.1}.
We denote by $Bj:K(\Z,3)\to BBU_\otimes$ a
first delooping of $j$, fitting in a homotopy commutative
diagram
\begin{equation}\label{kku-eq:sBj}
  \xymatrix{
  K(\Z,2)\ar[d]^{\ti s}\ar[r]^-{j} &
  BU_\otimes \ar[d]^{\ti s}\\
  \Omega K(\Z,3) \ar[r]^-{\Omega Bj} &
  \Omega BBU_\otimes \,,
  }
\end{equation}
where $\ti s$ denotes the homotopy equivalence 
which is right adjoint to the suspension
$s$ as in~\eqref{kku-eq:1squelet}.
We name $y_1\in\pi_2K(\Z,2)\cong\Z$ the generator that maps to
$y\in H_2(K(\Z,2);\F_p)$ by the Hurewicz homomorphism.
We have maps of based spaces 
\begin{equation*}
K(\Z,2)\stackrel{j}{\lr}
BU_\otimes\stackrel{c_0}{\lr} BU\times\{0\}\subset BU\times\Z \,,
\end{equation*}
where ${c_0}$ is the inclusion in $BU\times\Z$ followed by the translation
of the component of $1$ to that of $0$ in the $H$-group
$BU\times\Z$. The map $c_0j$
is a $\pi_2$-isomorphism, and we
define
\begin{equation*}
  u={c_0}_*j_*(y_1)\in\pi_2(BU\times\Z)\,. 
\end{equation*}
We call $u$ the {\sl Bott class}. We have an isomorphism of rings
\begin{equation*}
\pi_*(BU\times \Z)=\pi_*ku\cong\Z[u]
\end{equation*}
given by Bott periodicity.
The map ${c_0}_*:\pi_*(BU_\otimes)\to
\pi_*(BU\times \Z)$ is an isomorphism in positive degrees,
and we define $y_n\in \pi_{2n}(BU_\otimes)$ by requiring
${c_0}_*(y_n)=u^n$. 
Finally, we define 
\begin{equation}
\label{kku-eq:sigmaprime}
\sigma'_n\in V(1)_{2n+1}BBU_\otimes
\end{equation}
as the image of $y_n$ under the composition
\begin{equation*}
  \pi_{2n}BU_\otimes\stackrel{h_*}{\lr}
V(1)_{2n}BU_\otimes\stackrel{\sigma}{\lr}
V(1)_{2n+1}BBU_\otimes\,.
\end{equation*}
Here the first map is the Hurewicz homomorphism from
(unstable) homotopy to $V(1)$-homotopy, and
$\sigma$ is the suspension
induced by the map $s:\Sigma BU_\otimes\to
BBU_\otimes$.
\begin{lemma}\label{kku-lem:sigmas}
Consider the homomorphism
\begin{equation*}
Bj_*:V(1)_3K(\Z,3)\to V(1)_3 BBU_\otimes 
\end{equation*}
induced by the map defined above.
We have $\sigma'_1= (Bj)_*(e'_0)$, where 
$e'_0=\sigma^3(1)\in V(1)_{3}K(\Z,3)$, as given
in Definition~\textup{\ref{kku-def:prebott}}.
\end{lemma}
\begin{proof}
We have a commutative diagram 
\begin{equation*}
\xymatrix{
\pi_2K(\Z,2)\ar[r]^-{h_*}\ar[d]^{j_*} &
V(1)_2K(\Z,2)\ar[r]^-\sigma \ar[d]^{j_*} &
V(1)_3K(\Z,3)\ar[d]^{Bj_*} \\
\pi_2BU_\otimes \ar[r]^-{h_*} &
V(1)_2BU_\otimes \ar[r]^-\sigma  &
V(1)_3BBU_\otimes\,.
}
\end{equation*}
The right-hand square is induced in $V(1)$-homotopy from the
square left adjoint to the square~\eqref{kku-eq:sBj}.
The class $y_1\in  \pi_2K(\Z,2)$ was chosen so that
$h_*(y_1)=\sigma^2(1)$ in $V(1)_2K(\Z,2)
\cong H_2(K(\Z,2);\F_p)$. The lemma follows, since
\begin{equation*}
  \sigma'_1=\sigma h_*j_*(y_1)=(Bj)_*\sigma h_*(y_1) =
  (Bj)_*\sigma^3(1) =(Bj)_*(e'_0)\,.
\end{equation*}
\end{proof}
\par
The space $\Omega^\infty K(ku)$ is defined as the group completion of the
topological monoid $\coprod_{n}BGL_n(ku)$, with product modelling the
block-sum of matrices, see for instance \cite{Rkku-EKMM}*{VI.7}. 
The composition 
\begin{equation}\label{kku-eq:unit1}
w:BBU_\otimes\to BGL_1(ku)\to \coprod_{n}BGL_n(ku) \to \Omega^\infty K(ku)
\end{equation}
factors through an infinite
loop map $BBU_\otimes\to SL_1 K(ku)$, which is right adjoint to a map
\begin{equation*}
\omega:\Sigma^\infty_+BBU_\otimes\to K(ku)
\end{equation*}
of commutative $S$-algebras. 
We consider also the map of commutative $S$-algebras
\begin{equation*}
\phi:\Sigma^\infty_+K(\Z,3)\to
K(ku)
\end{equation*}
defined as the composition of the suspension of
$Bj:K(\Z,3)\to BBU_\otimes$ with the map~$\omega$.
\begin{definition}\label{kku-def:bott}
For $n\geq1$, we define
\begin{equation*}
\sigma_n=\omega_*(\sigma'_n)\in V(1)_{2n+1}K(ku)\,,
\end{equation*} 
where $\sigma'_n$ is the class given in~\eqref{kku-eq:sigmaprime}.
We define the ``higher Bott element'' as
\begin{equation*}
b=\phi_*(b')\in V(1)_{2p+2}K(ku)\,,
\end{equation*} 
where $b'\in V(1)_{2p+2} K(\Z,3)$ is the class given 
in Definition~\textup{\ref{kku-def:prebott}}.
\end{definition}
\begin{remark}
Notice that by Proposition~\ref{kku-prop:prebott} the classes $b$ and
$\sigma_1$ are related by a primary mod $(v_1)$ homotopy
Bockstein $\beta_{1,1}(b)=\sigma_1$.
\end{remark}
\begin{remark}\label{kku-rem:bott}
Assume that $p$ is an odd prime.
If $R$ is a number ring containing a primitive $p$-th root
of unity $\zeta_p$, for example $R=\Z[\zeta_p]$, then 
the mod $(p)$ algebraic $K$-theory of $R$ contains a
non-nilpotent class
\begin{equation*}
\beta\in V(0)_2 K(R)\,,
\end{equation*}
called the Bott element, which we referred to in the
introduction.
It was defined by Browder~\cite{Rkku-Br} using the composition
\begin{equation*}
  BC_p\to BGL_1R\to \Omega^\infty K(R)
\end{equation*}
analogous to~\eqref{kku-eq:unit1}, and its adjoint
\begin{equation*}
\phi:\Sigma^\infty_+ BC_p\to K(R)\,.
\end{equation*}
Here $C_p$ denotes the
cyclic subgroup of order $p$ of $GL_1(R)$ generated by
$\zeta_p$.
By inspection, the class $x=\zeta_p-1$ 
satisfies $x^p=0$ in the group-ring $\F_p[C_p]=V(0)_0C_p$, and has a
well defined ``transpotence'' $\beta'\in V(0)_2BC_p$,
supporting a primary mod $(p)$ homotopy Bockstein 
$\beta_{0,1}(\beta')\doteq \sigma(1)\in V(0)_1BC_p$. The classical Bott element 
can then be defined as 
\begin{equation*}
\beta= \phi_*(\beta')\in V(0)_2K(R)\,.
\end{equation*}
An embedding of rings $R\subset\C^\mathrm{top}$, where
$\C^\mathrm{top}$ has the Euclidean
topology, induces a map of commutative $S$-algebras
$\iota:K(R)\to K(\C^\mathrm{top})=ku$ 
in algebraic $K$-theory.
Browder's Proposition~\cite{Rkku-Br}*{2.2} implies that 
$\iota_*\phi_*(\beta')=u$, where $u$ is the Bott class
in $V(0)_*ku\cong P(u)$. This proves that $\beta$ is
non-nilpotent and is related to the Bott periodicity of
topological $K$-theory.  
Snaith showed~\cite{Rkku-Sn} that the relation 
$\beta'^p=v_1\beta'$ in $V(0)_*BC_p$ promotes
to the relation 
\begin{equation*}
\beta^{p-1}=v_1
\end{equation*}
in $V(0)_*K(R)$.
\end{remark}
\par
The remark above makes it clear that our construction of
$b\in V(1)_{2p+2}K(ku)$ is inspired from the
classical Bott element, and that these classes share
interesting properties. This provides some justification for
calling~$b$ a {\sl  higher Bott element}. Here {\sl higher}
refers to the fact that $b$ lives one chromatic step higher than 
$\beta$, in the sense that it is defined only in algebraic 
$K$-theory modulo $(p,v_1)$ and that it is related to
$v_2$-periodicity. Indeed, recall from
Theorem~\ref{kku-thm:main} and Proposition~\ref{kku-prop:integral} that $b$ 
is non-nilpotent and that the relation $b'^p=-v_2b'$ in
$V(1)_*K(\Z,3)$ promotes to the relation
\begin{equation*}
b^{p-1}=-v_2
\end{equation*}
in $V(1)_*K(ku)$.
Our proof of these assertions
relies on the computation of the cyclotomic trace for $ku$, 
and is much more technical then in the number ring
case\,: unfortunately, in the present situation we don't 
have an analogue of the map $K(R)\to
K(\C^\mathrm{top})$, but 
see the remark below for a possible candidate. 
\begin{remark}\label{kku-rem:algclosed}
John Rognes conjectured~\cite{Rkku-ARg} that if $\Omega_1$ is a separably 
closed $K(1)$-local pro-Galois extension of $ku$, 
in the sense of~\cite{Rkku-Ro08}, then
there is a weak equivalence
\begin{equation*}
L_{K(2)}K(\Omega_1)\simeq E_2\,,
\end{equation*}
where $L_{K(2)}$ is the Bousfield localization functor with
respect to the Morava $K$-theory $K(2)$, and
where $E_2$ is the second Morava $E$-theory spectrum~\cite{Rkku-GH} with
coefficients 
\begin{equation*}
(E_{2})_*=W(\F_{p^{2}})[[u_1]][u, u^{-1}]\,.
\end{equation*}
This would provide a map 
\begin{equation*}
\iota:K(ku)\to L_{K(2)}K(\Omega_1)\simeq E_2
\end{equation*}
that might play the role, at this chromatic level, of the
map $K(R)\to K(\C^{\mathrm{top}})$ mentioned in
Remark~\ref{kku-rem:bott}.
Since $V(1)_*E_2\cong \F_{p^2}[u, u^{-1}]$ with
$u^{p^2-1}=v_2$, we presume that the class $b$ would be
detected by the non-nilpotent class
\begin{equation*}
\iota_*(b)=\alpha u^{p+1} \in V(1)_*E_2
\end{equation*} 
for some $\alpha\in \F_{p^2}\setminus \F_p$ with
$\alpha^{p-1}=-1$. 
More generally, we expect that a periodic higher Bott element   
can be defined in $V(1)_*K(A)$ if $A$ is an commutative $S$-algebra
with an $S$-algebra map $A\to \Omega_1$ 
and a suitable $(p-1)$th-root of $v_1$ in $V(0)_*A$.
\end{remark}
\section{The trace map}
\label{kku-sec:trace}
In this section, we consider 
the B\"okstedt trace map~\cite{Rkku-BHM}
\begin{equation*}
\mathrm{tr}:K(ku) \to \thh(ku)
\end{equation*}
to topological Hochschild homology. This is a map of
commutative $S$-algebras, and it induces a homomorphism of
graded-commutative algebras in $V(1)$-homotopy, which we
just call the trace.
Our aim here is to prove that for $n\leq p-2$ the classes
$\sigma_n$ and $b$ defined above are non-zero in
$V(1)_*K(ku)$,
as well as some of their products, 
see Proposition~\ref{kku-prop:K-classes}.
We achieve this by showing that these classes have a non-zero trace in  
$V(1)_*\thh(ku)$.
To this end, we briefly recall the computation of
$V(1)_*\thh(ku)$ given in~\cite{Rkku-Au05}*{9.15}.
\par
The topological Hochschild homology  spectrum $\thh(ku)$ is a 
$ku$-algebra, and its $V(1)$-homotopy groups
form an algebra over the truncated polynomial algebra
$V(1)_*ku=P_{p-1}(u)$, where we also denote by $u$ the mod $(p,v_1)$
reduction of the Bott class $u\in\pi_2ku$.
There is a free $\F_p$-sub-algebra
$E(\lambda_1)\otimes P(\mu)$ in $V(1)_*\thh(ku)$, and there
is an isomorphism of
$E(\lambda_1)\otimes P(\mu)\otimes P_{p-1}(u)$-modules
\begin{equation}\label{kku-eq:thh}
V(1)_*\thh(ku)\cong E(\lambda_1)\otimes P(\mu)\otimes Q_*\,,
\end{equation}
where $Q_*$ is the $P_{p-1}(u)$-module given by
\begin{equation*}
Q_*=P_{p-1}(u)\oplus P_{p-2}(u)\{a_0,b_1,a_1,b_2,\dots,a_{p-2},b_{p-1}\}
\oplus P_{p-1}(u)\{a_{p-1}\}\,.
\end{equation*}
The degree of these generators is given by $|\lambda_1|=2p-1$,
$|\mu|=2p^2$, $|a_i|=2pi+3$ and $|b_j|=2pj+2$.
The isomorphism~\eqref{kku-eq:thh} is an isomorphism of $P_{p-1}(u)$-algebras if
the product on the $P_{p-1}(u)$-module generators of
$Q_*$ is given by the relations
\begin{equation}\label{kku-eq:thh2}
\begin{cases}
b_ib_j=ub_{i+j}& i+j\leq p-1\,,\\
b_ib_j=ub_{i+j-p}\mu& i+j\geq p\,,\\
a_ib_j=ua_{i+j}& i+j\leq p-1\,,\\
a_ib_j=ua_{i+j-p}\mu& i+j\geq p\,,\\
a_ia_j=0    &0\leq i,j\leq p-1\,.\\
\end{cases}
\end{equation}
Here by convention $b_0=u$.
For example we have a product
\begin{equation*}
(u^ka_i)(u^l b_j)=u^{p-2}a_{p-1}
\end{equation*} 
if $k+l=p-3$ and $i+j=p-1$.
\begin{remark}
The class $\mu$ is called $\mu_2$ in~\cite{Rkku-Au05},
but we adopt here the notation of~\cite{Rkku-AR02}.
\end{remark}
\par
The classes
$u^{n-1}a_0\in V(1)_{2n+1}\thh(ku)$ for $1\leq n\leq  p-2$ are
constructed as follows. 
The circle action $S^1_+\wedge \thh(ku)\to\thh(ku)$
restricts
in the homotopy category to a  map
$d:\Sigma\thh(ku)\to\thh(ku)$, which
in any homology theory $E_*$ induces Connes'
operator
\begin{equation}\label{kku-eq:connes}
d:E_*\thh(ku)\to E_{*+1}\thh(ku)\,.
\end{equation}
We have an $S$-algebra map $l:ku\to \thh(ku)$ given by the
inclusion of zero-simplices.
Composing the induced map in $E_*$-homology with
$d$ yields a suspension homomorphism
\begin{equation*}
dl_*:E_*ku\to E_{*+1}\thh(ku)\,,
\end{equation*}
see~\cite{Rkku-MS1}*{3.2} (it
is often denoted $\sigma$).
For $1\leq n\leq p-2$, we define the class $u^{n-1}a_0$ as the
image
\begin{equation*}
u^{n-1}a_0= dl_*(u^n)
\end{equation*}
of $u^{n}\in V(1)_{*}ku$.
Mapping to homology, we can show that these classes are
non-zero.
By Lemma~\ref{kku-lem:approx}, the Hurewicz homomorphism 
\begin{equation*}
h_*:V(1)_*\thh(ku)\to H_*(\thh(ku);\F_p)
\end{equation*}
is an isomorphism in degrees $*\leq 2p-3$ (notice that
$\alpha_1=0$ in $V(1)_*\thh(ku)$ since $\thh(ku)$ is a
$ku$-algebra). 
Let $x=h_*(u)\in H_2(ku;\F_p)$ be the image of $u\in
V(1)_2ku$. We then have $h_*(u^{n-1}a_0)=dl_*(x^n)$ in
$H_{2n+1}(\thh(ku);\F_p)$, 
and this class represents the permanent cycle $1\otimes
x^n\in E^1_{1,2n}(ku)$ in the B\"okstedt spectral sequence
\begin{equation*}
\begin{aligned}
E^1_{s,*}(ku)&=H_{*}(ku;\F_p)^{\otimes(s+1)}\,,\\
E^2_{s,*}(ku)&=HH^{\F_p}_{s,*}\big(H_{*}(ku;\F_p)\big)
\Rightarrow H_{s+*}(\thh(ku);\F_p)\,.\\
\end{aligned}
\end{equation*}
This proves that the classes
$h_*(u^{n-1}a_0)$ are non-zero
for these values of $n$.
We refer to~\cite{Rkku-Au05}*{\S9} for more details.
\begin{lemma}\label{kku-lem:tr}
If $1\leq n \leq p-2$, the class $\sigma'_n$
of~\textup{\eqref{kku-eq:sigmaprime}} 
maps to the class $u^{n-1}a_0$ under the composition
\begin{equation*}
V(1)_*BBU_\otimes
\stackrel{\omega_*}{\lr} 
V(1)_*K(ku)
\stackrel{{\mathrm{tr}}_*}{\lr}
V(1)_*\thh(ku)\,.
\end{equation*}
\end{lemma}
\begin{proof}
As mentioned above, $h_*:V(1)_{2n+1}\thh(ku)\to H_{2n+1}(\thh(ku);\F_p)$ is an
isomorphism for $n\leq p-2$ and maps $u^{n-1}a_0$ to
$dl_*(x^n)$.
Thus, passing to homology and using the definition of $\sigma'_n$ 
in~\eqref{kku-eq:sigmaprime}, if suffices to prove
that the composition 
\begin{equation*}
  H_{2n}(BU_\otimes;\F_p) \stackrel{\sigma}{\lr}
  H_{2n+1}(BBU_\otimes;\F_p) \stackrel{\tr_*\omega_*}{\lr}
  H_{2n+1}(\thh(ku);\F_p)
\end{equation*}
maps $z_n=h_*(y_n)\in H_{2n}(BU_\otimes;\F_p)$ to $dl_*(x^n)$.
Here we also denoted by $h_*$ the Hure\-wicz homomorphism
$\pi_{2n}BU_\otimes\to H_{2n}(BU_\otimes;\F_p)$. 
First, we need some information on the trace map. We will
use the following commutative diagram of spaces
\begin{equation}
  \xymatrix{
  BBU_\otimes\ar[r]^-i\ar[d]^{w} &
  B^{\mathrm{cy}}BU_\otimes\ar[d]^{\tau} &
  BU_\otimes\ar[l]_-{l}\ar[d]^{c_1} \\
  \Omega^\infty K(ku)\ar[r]^-{\Omega^\infty\tr} &
  \Omega^\infty \thh(ku) &
  \Omega^\infty ku\ar[l]_-{l}\,,
  }
  \label{kku-eq:trace}
\end{equation}
which is assembled from~\cite{Rkku-Sch04}*{\S4}.
Here the space $B^{\mathrm{cy}}BU_\otimes$ is the realization
of the cyclic nerve of the topological monoid $BU_\otimes$
and, as $\Omega^\infty \thh(ku)$, is equipped with a
canonical $S^1$-action.
The map $\tau$ is the realization of a morphism of
cyclic spaces, and is therefore
$S^1$-equivariant. The maps $l$ are given by the
inclusion of $0$-simplices, while $c_1$ is the inclusion of
the component of $1$.
There is a homotopy fibration~\cite{Rkku-Sch04}*{Proposition~3.1}
\begin{equation}\label{kku-eq:loopfibration}
  BU_\otimes \stackrel{l}{\lr} B^{\mathrm{cy}}BU_\otimes
  \stackrel{p}{\lr} BBU_\otimes\,,
\end{equation}
and the map $p$ admits a section up to homotopy $i: BBU_\otimes \to
B^{\mathrm{cy}}BU_\otimes$. 
\par
Let $d$ be Connes' operator on
$H_*(B^{\mathrm{cy}}BU_\otimes;\F_p)$ and 
$H_*(\Omega^\infty \thh(ku);\F_p)$. It
commutes with $\tau_*:H_*(B^{\mathrm{cy}}BU_\otimes;\F_p)\to
H_*(\Omega^\infty \thh(ku);\F_p)$ since $\tau$ is
equivariant.   
In the next lemma, we prove that 
\begin{equation*}
  d{l}_*(z_n) = i_*\sigma(z_n)
\end{equation*}
holds in $H_{2n+1}(B^{\mathrm{cy}}BU_\otimes;\F_p)$. 
Using~\eqref{kku-eq:trace}, we deduce
\begin{equation*}
  (\Omega^\infty\tr)_*w_*\sigma(z_n)=\tau_*i_*\sigma(z_n)=\tau_*dl_*(z_n)
  = d\tau_*l_*(z_n)=dl_*{c_1}_*(z_n)\,.
\end{equation*}
Finally, composing with the stabilization map 
\begin{equation*}
  {\mathrm{st}}:H_*(\Omega^\infty \thh(ku);\F_p)\to H_*(\thh(ku);\F_p)
\end{equation*}
to spectrum homology, we obtain
\begin{equation*}
  \tr_*\omega_*\sigma(z_n)={\mathrm{st}}(\Omega^\infty\tr)_*w_*\sigma(z_n)=
  {\mathrm{st}}dl_*{c_1}_*(z_n)=dl_*(x^n)\,.
\end{equation*}
For the last equality, we used that stabilization
commutes with $dl_*$, and that ${\mathrm{st}}{c_1}_*(z_n)=x^n$ for 
$1\leq n\leq p-2$.
\end{proof}
\begin{lemma}
  The equality $d{l}_*(z_n) = i_*\sigma(z_n)$ holds in 
  $H_{2n+1}(B^{\mathrm{cy}}BU_\otimes;\F_p)$.
    \label{kku-lem:primitive}
\end{lemma}
\begin{proof}
We consider the homotopy fibration~\eqref{kku-eq:loopfibration}. Since
$H_*(BU_\otimes;\F_p)$ is concentrated in even degrees, the
map $p_*:H_*(B^{\mathrm{cy}}BU_\otimes;\F_p)\to
H_*(BBU_\otimes;\F_p)$ restricts to an isomorphism
\begin{equation*}
p_*:{\mathrm{Prim}}\big(H_{2n+1}(B^{\mathrm{cy}}BU_\otimes;\F_p)\big)
\to {\mathrm{Prim}}\big(H_{2n+1}(BBU_\otimes;\F_p)\big)
\end{equation*}
of the subgroups of primitive elements in degree $2n+1$, with the
restriction of $i_*$ as inverse.
The class $l_*(z_n)$ is spherical, hence primitive, and it
follows from $d(1)=0$ that $dl_*(z_n)$ is also primitive.
\par
Next, 
we consider the diagram
\begin{equation*}
 \xymatrix{
 S^1\times BU_\otimes\ar@<-1ex>[d]\ar[r]^-{1\times l} &
 S^1\times B^{\mathrm{cy}}BU_\otimes\ar[r]^-{\mu} &
 B^{\mathrm{cy}}BU_\otimes\ar[d]^p \\
 S^1\wedge BU_\otimes \ar[rr]^{s}&&
 BBU_\otimes\,,
 }
\end{equation*}
where $\mu$ denotes the $S^1$-action on
$B^{\mathrm{cy}}BU_\otimes$ and $s$ the suspension 
map~\eqref{kku-eq:1squelet}.
This diagram is commutative, as can be checked at simplicial
level by using the definition of~$\mu$, see for
example~\cite{Rkku-Lo98}*{7.1.9}. Therefore
$p_*d{l}_*(z_n) = \sigma(z_n)$, and since $dl_*(z_n)$ is
primitive, we have
\begin{equation*}
d{l}_*(z_n) =  i_*p_*d{l}_*(z_n) =  i_*\sigma(z_n)\,.
\end{equation*}
\end{proof}
\begin{lemma}\label{kku-lem:b->b1}
The class $b'$
maps to the class $b_1$
under the composition
\begin{equation*}
V(1)_*K(\Z,3)
\stackrel{\phi_*}{\lr}V(1)_* K(ku)
\stackrel{{\mathrm{tr}}_*}{\lr}
V(1)_*\thh(ku)\,.
\end{equation*}
\end{lemma}
\begin{proof}
We know from Lemma~\ref{kku-lem:sigmas} and Lemma~\ref{kku-lem:tr}
that $e'_0\in
V(1)_{3}K(\Z,3)$ maps to the class
$a_0$ in $V(1)_{3}\thh(ku)$. We have primary
mod $(v_1)$ homotopy Bockstein 
\begin{equation*}
\beta_{1,1}(b')=e'_0
\ \ \textup{and} \ \ 
\beta_{1,1}(b_1)=a_0
\end{equation*}
in $V(1)_{*}K(\Z,3)$ and 
$V(1)_{*}\thh(ku)$ respectively, see Proposition~\ref{kku-prop:prebott} 
and~\cite{Rkku-Au05}*{9.19}. Moreover
$V(1)_{2p+2}\thh(ku)=\F_p\{b_1\}$, so that $\beta_{1,1}$ is
injective on this group. The result follows,
since
\begin{equation*}
\beta_{1,1}\tr_*\phi_*(b')=\tr_*\phi_*\beta_{1,1}(b')=\tr_*\phi_*(e'_0)=a_0\,.
\end{equation*}
\end{proof}
\par
Let $\kappa:ku\to ku_p$ be the completion at $p$. It induces
the inclusion $\Z[u]\to\Z_p[u]$ of coefficients rings.
\begin{definition}
We also denote by 
\begin{equation*}
\sigma_n\in V(1)_{2n+1}K(ku_p)\ \textup{ and }\  b\in V(1)_{2p+2}K(ku_p)
\end{equation*}
the image under $\kappa_*:V(1)_*K(ku)\to V(1)_*K(ku_p)$ of
the classes $\sigma_n$ and $b$ defined in~\ref{kku-def:bott}.
\end{definition}
\begin{proposition}\label{kku-prop:K-classes}
The classes
\begin{equation*}
\begin{cases}
b^{k} &\textup{for }0\leq k\leq p-2,\textup{ and} \\ 
\sigma_nb^l &\textup{for } 1\leq n \leq p-2\textup{ and }0\leq l\leq p-2-n
\end{cases}
\end{equation*}
are non-zero in $V(1)_*K(ku)$ and in $V(1)_*K(ku_p)$.
\end{proposition}
\begin{proof}
For  $V(1)_*K(ku)$, it follows from Lemma~\ref{kku-lem:tr}, 
Lemma~\ref{kku-lem:b->b1} and the
structure of $V(1)_*\thh(ku)$ given
in~\eqref{kku-eq:thh2}. In more detail, we have  
${\mathrm{tr}}_*(b^{k})=b_1^k\neq0$ for
$k\leq p-2$ and
${\mathrm{tr}}_*(\sigma_nb^{l})=u^{n-1}a_0b_1^l=u^{n+l-1}a_l\neq0$
for $l\leq p-3$ and $n+l-1\leq p-3$.
Notice that we have a commutative diagram
\begin{equation*}
\xymatrix{
V(1)_*K(ku)\ar[r]^-{{\mathrm{tr}}_*}\ar[d]_{\kappa_*}&
V(1)_*\thh(ku)\ar[d]^{\kappa_*}\\
V(1)_*K(ku_p)\ar[r]^-{{\mathrm{tr}}_*}&V(1)_*\thh(ku_p)\,.
}
\end{equation*}
The map $\kappa:\thh(ku)\to\thh(ku_p)$ is a weak equivalence
after $p$-completion, so in this diagram the right-hand $\kappa_*$ is 
an isomorphism. This
proves that the result also holds for $V(1)_*K(ku_p)$. 
\end{proof}
\begin{remark}
We claimed in Theorem~\ref{kku-thm:main} and
Proposition~\ref{kku-prop:integral} that $b$ is non-nilpotent in
$V(1)_*K(ku)$. However, we have
\begin{equation*}
  \mathrm{tr}_*(b^{p-1}) = \mathrm{tr}_*(b)^{p-1}=b_1^{p-1}=u^{p-2}b_{p-1}=0
\end{equation*}
in $V(1)_*\thh(ku)$, so that the B\"okstedt trace is not
sufficient for proving this assertion.
This is of course also predicted by our other claim that
$b^{p-1}=-v_2$ holds in
$V(1)_*K(ku)$. Indeed, $v_2$ maps to zero in
$V(1)_*\thh(ku)$ since $\thh(ku)$ is a $ku$-algebra. 
\end{remark}
\section{Algebraic $K$-theory in low degrees}
\label{kku-sec:low}
In this section, we compute the groups $V(1)_*K(ku_p)$
in degrees $*\leq 2p-2$. This complements the
computations presented in the next sections,  
which are based on evaluating the fixed points of 
$\thh(ku)$ and which are valid only in degrees larger than
$2p-2$, see Proposition~\ref{kku-prop:gamma1}.
\par
Consider the Adams summand 
\begin{equation*}
\ell_p=ku_p^{h\Delta}
\end{equation*} 
of $ku_p$,
where $\Delta\cong\Z/(p-1)$ is the finite subgroup of
the $p$-adic units, acting on $ku_p$ by $p$-adic Adams
operations, and where $(-)^{h\Delta}$ denotes the homotopy
fixed points.
By Theorem~10.2 of~\cite{Rkku-Au05}, the natural map 
$V(1)_*K(\ell_p)\to  V(1)_*K(ku_p)$ 
factors through an isomorphism 
\begin{equation}\label{kku-eq:delta-fix}
V(1)_*K(\ell_p)\cong \big(V(1)_*K(ku_p)\big)^\Delta\subset
V(1)_*K(ku_p)
\end{equation}
onto the elements of $V(1)_*K(ku_p)$ fixed under the induced
action of $\Delta$. In the sequel, we 
identify $V(1)_*K(\ell_p)$ with its image
in $V(1)_*K(ku_p)$.
\par
The $V(1)$-homotopy of $K(\ell_p)$ is computed
in~\cite{Rkku-AR02}.
In the degrees we are concerned with here, namely $*\leq
2p-2$, 
$V(1)_*K(\ell_p)$ is generated as an $\F_p$-vector space by
the classes listed in 
\begin{equation}\label{kku-eq:some-cl}
\{1,\lambda_1t^d,
s,\partial\lambda_1\,|\,0<d<p\}\,,
\end{equation}
of degree $|\lambda_1t^d|=2p-2d-1$, $|s|=2p-3$ and  
$|\partial\lambda_1|=2p-2$, see~\cite{Rkku-AR02}*{9.1} (where the
sporadic $v_2$-torsion class $s$ was denoted $a$).
The zeroth Postnikov section $\ell_p\to H\Z_p$ is a
$(2p-2)$-connected map, so that the induced map $K(\ell_p)\to
K(\Z_p)$ is $(2p-1)$-connected~\cite{Rkku-BoM94}*{Proposition~10.9}. 
All the classes listed in~\eqref{kku-eq:some-cl}
map to classes with same name in $V(1)_*K(\Z_p)$, which is
given by the formula
\begin{equation*}
V(1)_*K(\Z_p)\cong
E(\lambda_1)\oplus
\F_p\{s,\partial\lambda_1\}\oplus
\F_p\{\lambda_1t^d\,|\,0<d<p\}\,.
\end{equation*}
The name of the classes in this formula refers to
permanent cycles in the $S^1$ homotopy fixed-point spectral
sequence used in the computation of $V(1)_*K(\Z_p)$ by
traces, compare with Theorem~\ref{kku-thm:tcku}.
If desired, these classes could be given a more memorable name by
means of the inclusion
\begin{equation*}
V(0)_*K(\Z_p)\to V(0)_*K\big(\Q_p(\zeta_p)\big)\,,
\end{equation*}
in the target of which they can be decomposed as a product
of a unit and a
power of the Bott element $\beta\in
V(0)_2K\big(\Q_p(\zeta_p)\big)$.  
\par
Using the inclusion given in~\eqref{kku-eq:delta-fix}, we view the
classes listed in~\eqref{kku-eq:some-cl} as elements of
$V(1)_*K(ku_p)$. 
The following lemma implies that these classes are linearly
independent of the classes in $V(1)_*K(ku_p)$ constructed in the previous
section.
\begin{lemma}\label{kku-lem:non-fixed}
The non-zero classes $b^k$ and $\sigma_nb^l$ in
$V(1)_*K(ku_p)$ given in Pro\-po\-si\-tion~\textup{\ref{kku-prop:K-classes}} 
are not fixed under the action of $\Delta$.
\end{lemma}
\begin{proof}
All these classes map into $V(1)_*\thh(ku)$ to
classes which do not lie in the image of  $V(1)_*\thh(\ell_p)$, and hence
which are not fixed under the action of $\Delta$,
see Proposition~10.1 of~\cite{Rkku-Au05}.
\end{proof}
\begin{proposition}\label{kku-prop:low-deg}
The inclusion
\begin{equation*}
\F_p\{1,\sigma_n,\lambda_1t^d,s,\partial\lambda_1\,|\, 1\leq
n\leq p-2,\ 0<d<p\}\subset V(1)_*K(ku_p)
\end{equation*}
of graded $\F_p$-vector spaces
is an isomorphism in degrees $\leq 2p-2$.
\end{proposition}
\begin{proof}
We have constructed all the classes listed
above and have argued that they are linearly independent.
It suffices therefore to compute the dimension of
$V(1)_nK(ku_p)$ as an $\F_p$-vector space for all $0\leq n\leq2p-2$.
\par
Consider a double loop map $\Omega S^3\to BU_\otimes$
such that the composition 
\begin{equation*}
S^2\to \Omega S^3{\to}BU_\otimes\,,
\end{equation*}
where $S^2\to \Omega S^3$ is the adjunction unit,
represents the class $y_1\in\pi_2BU_\otimes$ defined in
Section~\ref{kku-sec:units}. By adjunction we have a map of $E_2$-ring spectra
\begin{equation*}
S[\Omega S^3]\to ku\,,
\end{equation*}
where $S[\Omega S^3]$ is another notation for the suspension
spectrum $\Sigma^\infty_+\Omega S^3$. We refer
to~\cite{Rkku-ARQ}*{Proposition~2.2} for some more details on the
construction of this map.
After $p$-com\-ple\-tion this map is $(2p-3)$-connected,
and induces a $(2p-2)$-connected map 
$K(S[\Omega S^3]_p)\to K(ku_p)$.
The dimension of the $\F_p$-vector space $V(1)_nK(S[\Omega S^3]_p)$ for
$n\leq 2p-2$ is computed in the following lemma, and this completes the proof of this
proposition. Notice that a priory 
\begin{equation*}
V(1)_{2p-2}K(S[\Omega S^3]_p)\to V(1)_{2p-2} K(ku_p)
\end{equation*}
is only surjective, but
luckily $V(1)_{2p-2}K(S[\Omega S^3]_p)$ is of rank one. 
Since we know that the rank of $V(1)_{2p-2} K(ku_p)$ is at least one, 
we also have an isomorphism in this degree.
\end{proof}
\begin{lemma}\label{kku-lem:KS3}
The dimension of $V(1)_nK(S[\Omega S^3]_p)$ as an $\F_p$-vector space is
\begin{equation*}
\begin{cases}
1&\textup{if }\ n=0, 1, 2p-2\,,  \\
2&\textup{if $n$ is odd with }\ 3\leq n \leq 2p-5\,, \\
3&\textup{if }\ n=2p-3\,,  \\
0&\textup{for other values of $n\leq 2p-2$\,.} 
\end{cases}
\end{equation*}
\end{lemma}
\begin{proof}
We compute $V(1)_*K(S[\Omega S^3]_p)$ in degrees less than
$2p-1$ by
using the cyclotomic trace map to topological
cyclic homology~\cite{Rkku-BHM}, which sits in a
cofibre sequence~\cite{Rkku-HM1}
\begin{equation*}
K(S[\Omega S^3]_p)_p\stackrel{\trc}{\lr}TC(S[\Omega S^3]_p)\to
\Sigma^{-1}H\Z_p\to \Sigma K(S[\Omega S^3]_p)_p\,.
\end{equation*}
Here $TC(X)=TC(X;p)$ denotes the ($p$-completed) topological cyclic
homology spectrum of a spectrum $X$.
By inspection, it suffices to prove that we
have
\begin{equation}\label{kku-eq:TCS3}
{\mathrm{dim}}_{\F_p}V(1)_nTC(S[\Omega S^3]_p)=
\begin{cases}
1&\textup{if }\ n=-1,0, 1, 2p-2\,,  \\
2&\textup{if $n$ is odd with }\ 3\leq n \leq 2p-3\,, \\
0&\textup{for other values of $n\leq2p-2$.} 
\end{cases}
\end{equation}
Indeed, $V(1)_*\Sigma^{-1}H\Z_p$ consists of a copy of
$\F_p$ in degrees $-1$ and $2p-2$, and is zero in other
degrees. We have an isomorphism $V(1)_{-1}TC(S[\Omega S^3]_p)\to 
V(1)_{-1}\Sigma^{-1}H\Z_p$, 
and the sporadic class $s$ is in the image of
the connecting homomorphism 
\begin{equation*}
V(1)_{2p-2}\Sigma^{-1}H\Z_p\to V(1)_{2p-3}K(S[\Omega
S^3]_p)\,,
\end{equation*}
by naturality with respect to
$S[\Omega S^3]_p\to H\Z_p$, see for example~\cite{Rkku-AR02}*{Proof
of~9.1}.
\par
The reduced topological cyclic homology spectrum
$\widetilde{TC}(S[\Omega S^3]_p)$ is the
homotopy fibre of the map $c:TC(S[\Omega S^3]_p)\to TC(S_p)$ induced by
the map $S^3\to *$ to a one-point space. The maps
$c$ admits a splitting, and we have a decomposition
\begin{equation*}
TC(S[\Omega S^3]_p)\simeq TC(S_p)\vee \widetilde{TC}(S[\Omega S^3]_p)\,.
\end{equation*}
The spectrum $TC(S_p)$ decomposes as
\begin{equation*}
TC(S_p)\simeq S_p\vee \Sigma\C P^\infty_{-1}\,,
\end{equation*}
where $\C P^\infty_{-1}$ is the ($p$-completed) Thom spectrum of minus the
canonical line bundle on $\C P^\infty$, see~\cite{Rkku-MS00}. 
The homology of $\Sigma\C P^\infty_{-1}$ is given by
\begin{equation*}
H_*(\Sigma\C
P^\infty_{-1};\F_p)\cong\F_p\{x_{i}\,|\,i\geq-1\}
\end{equation*}
with $|x_i|=2i+1$. Moreover these classes can be chosen so
that the relations 
\begin{equation*}
(P^1)^*(x_{p-2})=x_{-1}\ \textup{ and }\ 
(P^1)^*(x_{p-1})=0
\end{equation*}
hold. It follows from Lemma~\ref{kku-lem:approx} that
we have an inclusion
\begin{equation*}
\F_p\{c_{i}\,|\,-1\leq
i\leq{p-3}\}\cup\F_p\{\alpha(x_0)\}\subset 
V(1)_*(\Sigma\C P^\infty_{-1})\,,
\end{equation*}
which is an isomorphism in degrees $*\leq 2p-2$, with
$h_*(c_i)=x_i$. These classes have degree $|c_i|=2i+1$ and
$|\alpha(x_0)|=2p-2$.
\par
By~\cite{Rkku-BCCGHM}*{3.9}, we have a decomposition
\begin{equation*}
\widetilde{TC}(S[\Omega S^3]_p)\simeq \Sigma^\infty S^3_p\vee
\widetilde{V}\,,
\end{equation*}
where $\widetilde{V}$ is the ($p$-completed) homotopy fiber of the composition
\begin{equation*}
\Sigma^\infty\Sigma(ES^1_+\wedge_{S^1}LS^3) 
\stackrel{\mathrm{trf}}{\lr}
\Sigma^\infty LS^3 \stackrel{\epsilon_1}{\lr}
\Sigma^\infty S^3\,.
\end{equation*}
Here ${\mathrm{trf}}$ is the dimension-shifting
$S^1$-transfer on the free loop
space $LS^3$ of $S^3$, and $\epsilon_1$ is the evaluation
at $1\in S^1$, see~\cite{Rkku-MS00}.
We consider the Serre spectral sequence
\begin{equation*}
E^2_{**}=H_*\big(BS^1;H_*(LS^3,\F_p)\big)\Rightarrow
H_*(ES^1\times_{S^1}LS^3;\F_p)\,.
\end{equation*}
We have isomorphisms 
\begin{equation*}
  \begin{split}
H_*(BS^1;\F_p)&\ \cong H_*(K(\Z,2);\F_p)=\Gamma(y)\
\textup{ and}\\ 
H_*(LS^3;\F_p)&\ \cong P(z)\otimes E(dz)\,. 
  \end{split}
\end{equation*}
Here $z\in H_2(\Omega S^3;\F_p)\subset
H_2(L S^3;\F_p)$ and $dz\in H_3(L S^3;\F_p)$ is the
suspension of $z$ associated to the circle action on
$LS^3$. In particular, we have a non-zero $d^2$-differential
\begin{equation*}
d^2(yz)=dz\,.
\end{equation*}
For degree reasons no further non-zero differential involves the
classes in total degree less than $2p$, and we have an
inclusion
\begin{equation*}
P_p(y)\oplus\F_p\{
z^{j}\,|\, 1\leq j\leq p-1\}\subset H_*(ES^1\times_{S^1}LS^3;\F_p)
\end{equation*}
which is an isomorphism in degrees less than $2p$. We deduce
that the inclusion
\begin{equation*}
\Sigma\F_p\{z^j\,|\, 1\leq j\leq p-1\}\subset
H_*(\Sigma^\infty\Sigma(ES^1_+\wedge_{S^1}LS^3);\F_p)
\end{equation*}
is an isomorphism in degrees less the $2p-1$. 
The homomorphism
\begin{equation*}
(\epsilon_1{\mathrm{trf}})_*:H_*(\Sigma^\infty\Sigma(ES^1_+\wedge_{S^1}LS^3);\F_p)\to
H_*(S^3;\F_p)=E(e)
\end{equation*}
maps $\Sigma z$ to a generator $e$ of $H_3(S^3;\F_p)$ 
since the restriction of $\mathrm{trf}$
to $\Sigma^\infty\Sigma(S^1_+\wedge_{S^1} LS^3)$ is induced by
the circle action. This implies that we have an inclusion
\begin{equation*}
\F_p\{e, \Sigma z^{j}\,|\,2\leq j\leq p-2\}\subset
H_*(\Sigma^\infty S^3_p\vee\widetilde{V};\F_p)
\cong H_*(\widetilde{TC}(S[\Omega S^3]_p);\F_p)
\end{equation*}
which is an isomorphism in degrees smaller than $2p-1$. 
By Lemma~\ref{kku-lem:approx} 
\begin{equation*}
\F_p\{e, \Sigma z^{j}\,|\,2\leq j\leq p-2\}\subset
V(1)_*\widetilde{TC}(S[\Omega S^3]_p)
\end{equation*}
is also an isomorphism in degrees less than $2p-1$. 
In summary, we have
\begin{equation*}
V(1)_*TC(S[\Omega S^3]_p)\cong V(1)_*\oplus V(1)_*\Sigma\C
P^\infty_{-1}\oplus V(1)_*\widetilde{TC}(S[\Omega S^3]_p)\,,
\end{equation*}
which is isomorphic to
\begin{equation*}
\F_p\{1,\alpha_1,c_i,\alpha(x_0),e,\Sigma z^j\,|\,-1\leq
i\leq p-3,\ 2\leq j\leq p-2\}
\end{equation*}
in degrees smaller than $2p-1$. 
This proves that formula~\eqref{kku-eq:TCS3} for the rank of
the $\F_p$-vec\-tor space $V(1)_*TC(S[\Omega S^3]_p)$ is correct.
\end{proof}
\begin{remark}
In an earlier proof of this lemma we used the space $BU(1)$
and the map $\theta:\Sigma_+^\infty BU(1)\to ku$ of
commutative $S$-algebras. 
I thank John Rognes for noticing that using $\Omega S^3$
instead simplifies the computation. The maps 
\begin{equation*}
S[\Omega S^3]\to \Sigma_+^\infty BU(1)\to ku
\end{equation*}
are $\pi_0$-isomorphisms and rational equivalences. We use this in~\cite{Rkku-ARQ} to
determine the rational homotopy type of $K(ku)$. 
\end{remark}
\section{The fixed points}
\label{kku-sec:fixed}
In this section we compute the $V(1)$-homotopy groups of the
homotopy limit
\begin{equation*}
TF(ku_p)=\holim_{n,F}\thh(ku_p)^{C_{p^n}}\,,
\end{equation*}
where $F:\thh(ku_p)^{C_{p^{n+1}}}\to \thh(ku_p)^{C_{p^n}}$ is
the Frobenius map.
This will be used in the next section to compute the
topological cyclic homology of $ku_p$.
The strategy to perform such computations was developed in
\citelist{\cite{Rkku-BoM94} \cite{Rkku-Ts} \cite{Rkku-HM1} }, but we will
closely follow the exposition and adopt the notations of
\cite{Rkku-AR02}*{\S3, \S5 and \S6}, with an exception\,: 
the $G$ Tate construction on an equivariant $G$ spectrum $X$ will be denoted
by $X^{tG}$ instead of $\hat{\mathbb H}(G,X)$. 
We refer the reader to \cite{Rkku-AR02}*{\S3} for a brief review of
the homotopy commutative norm-restriction diagram 
\begin{equation*}
\xymatrix@C=4ex@R=4ex{
 & K(ku_p) \ar[d]_{\tr_n} \ar[dr]^{\tr_{n-1}} \\
        \thh(ku_p)_{hC_{p^n}} \ar[r]^-N \ar@{=}[d] &
	\thh(ku_p)^{C_{p^n}} \ar[r]^-R \ar[d]^{\Gamma_n} &
	\thh(ku_p)^{C_{p^{n-1}}}\ar[r] \ar[d]^{\hat\Gamma_n} &
        \Sigma \thh(ku_p)_{hC_{p^n}}\ar@{=}[d] \\
\thh(ku_p)_{hC_{p^n}} \ar[r]^-{N^h} &
	\thh(ku_p)^{hC_{p^n}} \ar[r]^-{R^h} &
	\thh(ku_p)^{tC_{p^n}} \ar[r] &
        \Sigma \thh(ku_p)_{hC_{p^n}}
}
\end{equation*}
for any $n\geq1$, which is our essential tool. 
By passage to homotopy limits over the Frobenius maps, we
obtain the homotopy commutative diagram
\begin{equation*}
\xymatrix@C=4ex@R=4ex{
 & K(ku_p) \ar[d]_{\tr_F} \ar[dr]^{\tr_F} \\
        \Sigma\thh(ku_p)_{hS^1} \ar[r]^-N \ar@{=}[d] &
	TF(ku_p) \ar[r]^-R \ar[d]^{\Gamma} &
	TF(ku_p) \ar[r] \ar[d]^{\hat\Gamma} &
        \Sigma^2 \thh(ku_p)_{hS^1}\ar@{=}[d] \\
        \Sigma\thh(ku_p)_{hS^1} \ar[r]^-{N^h} &
	\thh(ku_p)^{hS^1} \ar[r]^-{R^h} &
	\thh(ku_p)^{tS^1} \ar[r] &
        \Sigma^2 \thh(ku_p)_{hS^1}\,.
}
\end{equation*}
The map $i_*:V(1)_*\thh(\ell_p)\to\thh(ku_p)$ factors through an
isomorphism onto the $\Delta$-fixed elements of
$V(1)_*\thh(ku_p)$,
\begin{equation}\label{kku-eq:thh-desc}
i_*:V(1)_*\thh(\ell_p)\stackrel{\cong}{\lr}\big(V(1)_*\thh(ku_p)\big)^\Delta
\subset V(1)_*\thh(ku_p)\,,
\end{equation}
see~\cite{Rkku-Au05}*{10.1}. 
The corresponding results hold also for
the $C_{p^n}$ or $S^1$ homotopy fixed points of
$\thh$, for the $C_{p^n}$ or $S^1$ Tate construction
on $\thh$, and for  $\tc$ and $K$,
see \cite{Rkku-Au05}*{10.2}. In the sequel, we identify
$V(1)_*\thh(\ell_p)$, $V(1)_*\tc(\ell_p)$, etc. with 
their image under $i_*$.
We have a similar statement for the various spectral
sequences computing the $V(1)$-homotopy of these spectra.
\begin{lemma}\label{kku-lem:ss-desc}
Let $G=S^1$ or $G=C_{p^n}$, and let $E^*(G,\ell_p)$ and
$E^*(G,ku_p)$ be the $G$ homotopy fixed-point spectral
sequences converging strongly to $V(1)_*\thh(\ell_p)^{hG}$ and
to $V(1)_*\thh(ku_p)^{hG}$, respectively.
Then the morphism of spectral sequences
induced by the map $\ell_p\to ku_p$ is equal to the inclusion
of the $\Delta$ fixed points
\begin{equation*}
E^*(G,\ell_p)=\big(E^*(G,ku_p)\big)^\Delta\subset
E^*(G,ku_p)\,.
\end{equation*}
This holds also for the morphism induced on the $G$ Tate spectral
sequences converging to $V(1)_*\thh(\ell_p)^{tG}$ and
$V(1)_*\thh(ku_p)^{tG}$, which is given
by
\begin{equation*}
\hat E^*(G,\ell_p)=\big(\hat E^*(G,ku_p)\big)^\Delta\subset
\hat E^*(G,ku_p)\,.
\end{equation*}
\end{lemma}
\begin{proof} The group $\Delta$ acts on $ku_p$ by
$S$-algebra maps, and it acts $S^1$-equivariently on
$\thh(ku_p)$. In particular $\Delta$ acts by morphisms of
spectral sequences on $E^*(G,ku_p)$ and $\hat E^*(G,ku_p)$, 
and hence it suffices to prove that the
claims hold at the level of the $E^2$-terms. This follows
from~\eqref{kku-eq:thh-desc}. 
\end{proof}
\par
From now on, we will omit $ku_p$ from the notation and just
write $E^*(G)$ and $\hat E^*(G)$ for the $G$ homotopy 
fixed-point and Tate spectral sequences
converging to $V(1)_*\thh(ku_p)^{hG}$ and
$V(1)_*\thh(ku_p)^{tG}$, respectively.
\par
At this point, we recall the notion of $\delta$-weight
introduced in~\cite{Rkku-Au05}*{8.2}. We fix a generator $\delta$
of the group $\Delta$ acting on $ku_p$, $K(ku_p)$,
$\thh(ku_p)$, $\tc(ku_p)$, etc. The self-map
$\delta_*$ of $V(1)_*ku_p=P_{p-1}(u)$ maps $u$ to $\alpha u$ for some
generator $\alpha$ of $\F_p^\times$. 
We say that a class
$v\in V(1)_*K(ku_p)$ has $\delta$-weight
$i\in\Z/(p-1)$ if $\delta_*(v)=\alpha^iv$. 
The same convention holds
for classes in $V(1)_*\thh(ku_p)$, $V(1)_*\tc(ku_p)$, etc.
For example, the generators $a_i$ and $b_j$ of $V(1)_*\thh(ku_p)$ given
in~\eqref{kku-eq:thh} all have $\delta$-weight $1$,
see~\cite{Rkku-Au05}*{10.1}. Similarly, it follows from its
definition that $b\in V(1)_*K(ku_p)$ has $\delta$-weight
$1$. Since $\delta_*$ is diagonalizable,  we can
reinterpret Lemma~\ref{kku-lem:ss-desc} by saying that each of these spectral
sequences for $ku_p$ has an extra $\Z/(p-1)$-grading given by the
$\delta$-weight, and that its homogeneous summand of
$\delta$-weight $0$ consists of the corresponding spectral sequence for
$\ell_p$. Together with the internal and filtration degrees,
the $\delta$-weight endows the $E^r$-terms of these spectral
sequences with a tri-grading that we will refer to in 
the computations below.
\par
By a computation of McClure and Staffeldt~\cite{Rkku-MS1}, \cite{Rkku-AR02}*{2.6},
we have an isomorphism of $\F_p$-algebras
\begin{equation*}
V(1)_*\thh(\ell_p)\cong E(\lambda_1,\lambda_2)\otimes
P(\mu)\,.
\end{equation*}
The induced map $V(1)_*\thh(\ell_p)\to
V(1)_*\thh(ku_p)$ sends $\lambda_1$ and
$\mu$ to the classes with same name, and $\lambda_2$ to the
class $a_1 b_1^{p-2}$. 
\begin{remark}
In the sequel, we will frequently denote by
$\lambda_2$ the class $a_1b_1^{p-2}$.
\end{remark}
\par
The $C_p$-Tate spectral sequence
\begin{equation*}
\hat E(C_p)_{s,t}^2=\hat
H^{-s}\big(C_p,V(1)_t\thh(ku_p)\big)\Rightarrow
V(1)_{s+t}\thh(ku_p)^{tC_p}
\end{equation*}
has an $E_2$-term given by
\begin{equation*}
\hat E(C_p)^2=P(t,t^{-1})\otimes E(u_1)\otimes V(1)_*\thh(ku_p)
\end{equation*}
with $t$ in bidegree $(-2,0)$, $u_1$ in bidegree $(-1,0)$, and
$w\in V(1)_t\thh(ku_p)$ in bidegree $(0,t)$.
Recall the description of $V(1)_*\thh(ku_p)$ given
in~\eqref{kku-eq:thh}.
\begin{lemma}\label{kku-lem:cp-tate}
In the $C_{p}$ Tate spectral sequence $\hat E^*(C_{p})$
the classes $\lambda_1$, $\lambda_2$, $b_1$ and $t\mu$ are infinite cycles.
There are non-zero differentials
\begin{equation*}
\begin{split}
d^2(b_i)&=(1-i)a_i t\\
d^{2p}(t^{1-p})&\doteq \lambda_1\cdot t\\
d^{2p^2}(t^{p-p^2})&\doteq \lambda_2\cdot t^p\\
d^{2p^2+1}(u_1\cdot t^{-p^{2}})&\doteq  t\mu
\end{split}
\end{equation*}
with $0\leq i\leq p-1$. The spectral sequence collapses at the
$\hat E^{2p^2+2}$-term, leaving
\begin{equation*}
\begin{split}
\hat E^{\infty}(C_{p})=\ &
P(t^{\pm p^2}) \otimes E(\lambda_1,a_1) \otimes P_{p-1}(b_1)
\\
& \oplus E(\lambda_1) \otimes P_{p-2}(b_1) \otimes 
\F_p\{a_1t^j,b_1t^j\,|\,v_p(j)=1\}\,.
\end{split}
\end{equation*}
\end{lemma}
\begin{remark}
Beware that in the lemma above, the index $j$ appearing as a 
power of $t$ runs over all integers,
positive or negative, with specified $p$-adic valuation. The
same remark holds for
the Lemmas~\ref{kku-lem:cpn-tate} and~\ref{kku-lem:tate-S1} below,
and also for the power $j$ of $\mu$ in Lemmas~\ref{kku-lem:cpn-fixed}
and~\ref{kku-lem:fixed-S1} below.
\end{remark}
\begin{proof}
We know from~\cite{Rkku-AR02}*{Proposition~4.8} that $t\mu$ is
an infinite cycle. The classes $\lambda_1$, $\lambda_2$ and
$b_1$ are also infinite cycles, see the argument given 
at the top of~\cite{Rkku-AR02}*{page~21}. 
\par
Let $d$ be Connes' operator~\eqref{kku-eq:connes} on
$V(1)_*\thh(ku_p)$, and
recall from above the notation $b_0=u$. We have
\begin{equation*}
d(b_0)=a_0\,,
\end{equation*}
and this relation is detected via the Hurewicz
homomorphism in mod $(p)$ homology, see~\cite{Rkku-Au05}*{\S9}.
It follows from~\cite{Rkku-Ro98}*{3.3} that 
in the $S^1$ homotopy fixed-point spectral sequence
\begin{equation*}
E^2(S^1)=P(t)\otimes V(1)_*\thh(ku_p)\Rightarrow
V(1)_*\thh(ku_p)^{hS^1}
\end{equation*}
we have a $d^2$-differential
\begin{equation*}
d^2(b_0)=a_0 t\,.
\end{equation*}
Since $E^2(S^1)$ injects into $\hat E^2(C_p)$ via $R^hF$, 
this differential is also present in $\hat
E^2(C_p)$. The differentials $d^2(b_i)=(1-i)a_i t$ for
$i\neq0$ follow easily from the case $i=0$ and the multiplicative
structure.  Indeed $d^2(\mu)=0$ for
degree reasons, and hence
$d^2(u^2 \mu)=2u \mu  a_0 t$. From
the relation $b_i b_{p-i}=u^2\mu$ we deduce that 
$d^2(b_i)=\alpha_ia_i t$ for some $\alpha_i\in\F_p$, 
because in $V(1)_*\thh(ku_p)$ the equation $x b_{p-i}=u\mu a_0$ 
has $x=a_i$ as unique (homogeneous) solution. 
First, notice that
$0=d^2(b_1^{p-1})=(p-1)\alpha_1\lambda_2$, so we have $\alpha_1=0$. 
Next, the relation $b_1 b_{p-1}=u^2 \mu$ implies that $\alpha_{p-1}=2$, while
$b_1 b_i=u b_{i+1}$ for $i\leq p-2$ implies that
$\alpha_i=1+\alpha_{i+1}$. We deduce that $\alpha_i=1-i$,
proving the claim on the $d^2$-differential,
which leaves
\begin{equation*}
\hat E^3(C_{p})=P(t^{\pm1},t\mu)\otimes E(u_1,\lambda_1,a_1)\otimes
P_{p-1}(b_1)\,.
\end{equation*}
Lemma~\ref{kku-lem:ss-desc} determines the given next three
non-zero differentials, by comparison with the case of the $\ell_p$
treated in~\cite{Rkku-AR02}*{5.5}, and this takes care of the
summand of $\delta$-weight zero. 
The only algebra generators of $\hat E^3(C_{p})$ 
of non-zero $\delta$-weight are $a_1$ and $b_1$. We know
that $b_1$ is an infinite cycle.
In the $S^1$ Tate spectral sequence, using the known
differentials, the tri-grading and the product, it is easy to see that
$a_1$ survives to the $E^{2p^2+2}$-term. Therefore $a_1$
also survives to the $E^{2p^2+2}$-term
in $\hat E^*(C_{p})$, via the morphism of
spectral sequences induced by $F$. The $d^{2p}$ differential
leaves 
\begin{equation*}
\hat E^{2p+1}(C_{p})=P(t^{\pm p},t\mu)\otimes E(u_1,\lambda_1,a_1)\otimes
P_{p-1}(b_1)\,,
\end{equation*}
and the $d^{2p^2}$ differential leaves
\begin{equation*}
\begin{split}
\hat E^{2p^2+1}(C_{p})=\ &
P(t^{\pm p^2},t\mu) \otimes E(u_1,\lambda_1,a_1) \otimes P_{p-1}(b_1)
\\
& \oplus E(u_1,\lambda_1)\otimes P_{p-2}(b_1) \otimes P(t\mu) \otimes 
\F_p\{a_1t^j,b_1t^j\,|\,v_p(j)=1\}\,,
\end{split}
\end{equation*}
as can be computed using the relation $a_1\cdot b^{p-2}=\lambda_2$.
Finally, $d^{2p^2+1}$
leaves
\begin{equation*}
\begin{split}
\hat E^{2p^2+2}(C_{p})=\ &
P(t^{\pm p^2}) \otimes E(\lambda_1,a_1) \otimes P_{p-1}(b_1)
\\
& \oplus E(\lambda_1) \otimes P_{p-2}(b_1) \otimes 
\F_p\{a_1t^j,b_1t^j\,|\,v_p(j)=1\}\,,
\end{split}
\end{equation*}
and at this stage the spectral sequence collapses for
bidegree reasons.
\end{proof}
\begin{remark}
The $d^2$-differential can also be determined by computing
$d(b_i)$ for $i\geq0$, using Connes' operator in 
Hochschild homology (c.f.~\cite{Rkku-Au05}*{3.4}). 
\end{remark}
\begin{definition}
We call a homomorphism of graded groups $k$-coconnected if it
is an isomorphism in all dimensions greater than $k$ and
injective in dimension $k$.
\end{definition}
\begin{proposition}\label{kku-prop:gamma1}
The algebra map 
\begin{equation*}
(\hat\Gamma_1)_*:V(1)_*\thh(ku_p)\to 
V(1)_*\thh(ku_p)^{tC_p}
\end{equation*} 
factorizes as the localization away
from $\mu$, followed by an isomorphism
\begin{equation*}
V(1)_*\thh(ku_p)[\mu^{-1}]\to V(1)_*\thh(ku_p)^{tC_p}
\end{equation*}
given by 
\begin{equation*}
\lambda_1\mapsto \lambda_1,\ \mu\mapsto t^{-p^2},\ 
b_i\mapsto t^{(1-i)p}b_1,\ \hbox{and}\ a_i\mapsto
t^{(1-i)p}a_1
\end{equation*} 
for $0\leq i\leq p-1$, up to some non-zero
scalar multiples. 
In particular the map $(\hat\Gamma_1)_*$  is $(2p-2)$-coconnected.
\end{proposition}
\begin{proof}
By naturality with respect to $\ell_p\to ku_p$ and by the computation of
$(\hat\Gamma_1)_*$ for $\ell_p$ given in~\cite{Rkku-AR02}*{Theorem~5.5}, we
know that the map $(\hat\Gamma_1)_*$ for $ku_p$ satisfies 
\begin{equation*}
\lambda_1\mapsto \lambda_1,\ \ \lambda_2\mapsto \lambda_2\
\hbox{ and }\ \mu\mapsto t^{-p^2}\,.
\end{equation*}
In $V(1)_*\thh(ku_p)$ we have multiplicative relations
$u^{p-3}a_i b_j=\lambda_2$ for $i+j=p-1$, from which we deduce that
$(\hat\Gamma_1)_*(u^ka_i)\neq0$ and
$(\hat\Gamma_1)_*(u^kb_i)\neq0$ for any $0\leq k\leq p-3$ and any
$0\leq i\leq p-1$. For degree reasons, this forces
\begin{equation*}
(\hat\Gamma_1)_*(a_i)=t^{(1-i)p}a_1\ \hbox{ and }\
(\hat\Gamma_1)_*(b_i)=t^{(1-i)p}b_1
\end{equation*}
up to some non-zero scalar multiples.
\end{proof}
\begin{corollary}\label{kku-cor:gammas}
The canonical maps 
\begin{equation*}
\begin{split}
\Gamma_n:\ &\thh(ku_p)^{C_{p^n}}\to\thh(ku_p)^{hC_{p^n}}\,,\\
\hat\Gamma_n:\
&\thh(ku_p)^{C_{p^{n-1}}}\to\thh(ku_p)^{tC_{p^n}}\,,\\
\Gamma:\ &TF(ku_p)\to\thh(ku_p)^{hS^1}\,,\\
\hat\Gamma:\ &TF(ku_p)\to\thh(ku_p)^{tS^1}\,,
\end{split}
\end{equation*}
for $n\geq1$ all induce $(2p-2)$-coconnected maps in
$V(1)$-homotopy.
\end{corollary}
\begin{proof}
The claims for $\Gamma_n$ and $\hat\Gamma_n$ follow from
Proposition~\ref{kku-prop:gamma1} and the generalization of a theorem of
Tsalidis~\cite{Rkku-Ts} given in~\cite{Rkku-BBLR}. The claims for $\Gamma$ and
$\hat\Gamma$ follow by passage to homotopy limits.
\end{proof}
\begin{definition}\label{kku-def:rn}
Let $r(n)=0$ for all $n\leq0$, and let $r(n)=p^n+r(n-2)$ for
all $n\geq1$. Thus $r(2n-1)=p^{2n-1}+\dots+p$ (odd powers)
and $r(2n)=p^{2n}+\dots+p^2$ (even powers).
\end{definition}
\begin{lemma}\label{kku-lem:cpn-tate}
In the $C_{p^n}$ Tate spectral sequence $\hat E^*(C_{p^n})$
the classes $\lambda_1$, $\lambda_2$, $b_1$ and $t\mu$ are infinite cycles.
There are non-zero differentials
\begin{equation*}
\begin{split}
d^2(b_i)&=(1-i)a_i t\\
d^{2p}(t^{1-p})&\doteq \lambda_1\cdot t\\
d^{2p^2}(t^{p-p^2})&\doteq \lambda_2\cdot t^p
\end{split}
\end{equation*}
with $0\leq i\leq p-1$, leaving
\begin{equation*}
\begin{split}
\hat E^{2p^2+1}(C_{p^n})=\ &
P(t^{\pm p^2}) \otimes E(u_n, \lambda_1,a_1) \otimes P_{p-1}(b_1)
\otimes P(t\mu)\\
& \oplus E(u_n, \lambda_1) \otimes P_{p-2}(b_1) \otimes P(t\mu)
\otimes \F_p\{a_1t^j,b_1t^j\,|\,v_p(j)=1\}\,.
\end{split}
\end{equation*}
If $n\geq2$, then for each $1\leq k\leq n-1$ there is a
triple of non-zero differentials
\begin{equation*}
\begin{split}
d^{2r(2k)+2}(b_1t^j) & \doteq 
a_1t^j \cdot t^{p^{2k}} \cdot (t\mu)^{r(2k-2)+1}\\
d^{2r(2k+1)}(t^{p^{2k}-p^{2k+1}}) & \doteq  
\lambda_1 \cdot t^{p^{2k}} \cdot (t\mu)^{r(2k-1)}\\
d^{2r(2k+2)}(t^{p^{2k+1}-p^{2k+2}}) & \doteq  
\lambda_2 \cdot t^{p^{2k+1}} \cdot (t\mu)^{r(2k)}\\
\end{split}
\end{equation*}
with $v_p(j)=2k-1$, leaving
\begin{equation*}
\begin{split}
\hat E^{2r(2k+2)+1}(C_{p^n})=\ &
P(t^{\pm p^{2k+2}}) \otimes E(u_n, \lambda_1,a_1) \otimes P_{p-1}(b_1)
\otimes P(t\mu)\\
& \oplus E(u_n, \lambda_1) \otimes P_{p-2}(b_1) \otimes P(t\mu)
\otimes \F_p\{a_1t^j,b_1t^j\,|\,v_p(j)=2k+1\}\\
& \oplus \bigoplus_{1\leq m\leq k} \hat T_m(C_{p^n})\,,
\end{split}
\end{equation*}
where
\begin{equation*}
\begin{split}
\hat T_m(C_{p^n})=\ &
E(u_n,\lambda_1) \otimes P_{r(2m)}(t\mu)\otimes
\F_p\{\lambda_2t^j\,|\,v_p(j)=2m+1\}\\
&\oplus E(u_n,a_1) \otimes P_{p-1}(b_1) \otimes P_{r(2m-1)}(t\mu)\otimes
\F_p\{\lambda_1t^j\,|\,v_p(j)=2m\}\\
&\oplus E(u_n,\lambda_1) \otimes P_{p-2}(b_1) \otimes P_{r(2m-2)+1}(t\mu)\otimes
\F_p\{a_1t^j\,|\,v_p(j)=2m-1\}\,.
\end{split}
\end{equation*} 
For $n\geq1$, there is a last non-zero differential
\begin{equation*}
d^{2r(2n)+1}(u_n\cdot t^{-p^{2n}})  \doteq (t\mu)^{r(2n-2)+1}
\end{equation*}
after which the spectral sequence collapses, leaving
\begin{equation*}
\begin{split}
\hat E^{\infty}(C_{p^n})=\ &
P(t^{\pm p^{2n}}) \otimes E(\lambda_1,a_1) \otimes P_{p-1}(b_1)
\otimes P_{r(2n-2)+1}(t\mu)\\
&\oplus E(\lambda_1) \otimes P_{p-2}(b_1) \otimes P_{r(2n-2)+1}(t\mu)
\otimes
\F_p\{a_1t^j,b_1t^j\,|\,v_p(j)=2n-1\}\\
& \oplus \bigoplus_{1\leq m\leq n-1} \hat T_m(C_{p^n})\,.
\end{split}
\end{equation*} 
\end{lemma}
\par
Next, we describe the $C_{p^n}$ homotopy
fixed-point spectral sequence $E^*(C_{p^n})$.
It is algebraically easier to describe the $E^r$-terms of
the $C_{p^n}$ homotopy fixed-point spectral
sequence for $\thh(ku_p)^{tC_p}$, which we denote abusively
by
\begin{equation*}
 \mu^{-1}E^*(C_{p^n})\Rightarrow
 V(1)_*(\thh(ku_p)^{tC_p})^{hC_{p^n}}\,,
\end{equation*}
compare with~\cite{Rkku-AR02}*{page~23}.
We know from Proposition~\ref{kku-prop:gamma1} that the
map 
\begin{equation*}
\hat\Gamma_1^{hC_{p^n}}:\thh(ku_p)^{hC_{p^n}}\to
(\thh(ku_p)^{tC_p})^{hC_{p^n}}
\end{equation*}
induces a
morphism of spectral sequences
\begin{equation*}
  E^*(C_{p^n})\to \mu^{-1}E^*(C_{p^n})
\end{equation*}
which on $E^2$-terms (but not on higher terms)
indeed corresponds to inverting $\mu$.
By the same Proposition and by strong convergence of the spectral
sequences, the map $\hat\Gamma_1^{hC_{p^n}}$ induces a
$(2p-2)$-coconnected homomorphism in $V(1)$-homotopy. 
\begin{lemma}\label{kku-lem:cpn-fixed}
In the $C_{p^n}$ homotopy fixed-point spectral sequence $\mu^{-1}E^*(C_{p^n})$
the classes $\lambda_1$, $\lambda_2$, $b_1$ and $t\mu$ are infinite cycles.
There are non-zero differentials
\begin{equation*}
\begin{split}
d^2(b_i)&=(1-i)a_i t\\
d^{2p}(\mu^{p-1})&\doteq \lambda_1\cdot\mu^{-1}\cdot(t\mu)^p\\
d^{2p^2}(\mu^{p^2-p})&\doteq
\lambda_2\cdot\mu^{-p}\cdot(t\mu)^{p^2}
\end{split}
\end{equation*}
with $0\leq i\leq p-1$, leaving
\begin{equation*}
\begin{split}
\mu^{-1}E^{2p^2+1}(C_{p^n})=\ &
P(\mu^{\pm p^2}) \otimes E(u_n, \lambda_1,a_1) \otimes P_{p-1}(b_1)
\otimes P(t\mu)\\
& \oplus E(u_n, \lambda_1) \otimes P_{p-2}(b_1) \otimes P(t\mu)
\otimes \F_p\{a_1\mu^j,b_1\mu^j\,|\,v_p(j)=1\}\\
& \oplus T_1(C_{p^n})\,,
\end{split}
\end{equation*}
where
\begin{equation*}
\begin{split}
T_1(C_{p^n})=\ &
E(u_n,\lambda_1) \otimes P_{p^2}(t\mu)\otimes
\F_p\{\lambda_2\mu^j\,|\,v_p(j)=1\}\\
&\oplus E(u_n,a_1) \otimes P_{p-1}(b_1) \otimes P_{p}(t\mu)\otimes
\F_p\{\lambda_1\mu^j\,|\,v_p(j)=0\}\\
&\oplus E(u_n,\lambda_1) \otimes P_{p-2}(b_1) \otimes P(\mu^{\pm1})\otimes
\F_p\{a_i\,|\,0\leq i\leq p-1,\ i\neq1\}\,.
\end{split}
\end{equation*}
If $n\geq2$, then for each $2\leq k\leq n$ there is a triple
of non-zero differentials
\begin{equation*}
\begin{split}
d^{2r(2k-2)+2}(b_1\mu^j) & \doteq 
a_1\mu^j \cdot\mu^{-p^{2k-2}} \cdot (t\mu)^{r(2k-2)+1}\\
d^{2r(2k-1)}(\mu^{p^{2k-1}-p^{2k-2}}) & \doteq  
\lambda_1 \cdot\mu^{-p^{2k-2}} \cdot (t\mu)^{r(2k-1)}\\
d^{2r(2k)}(\mu^{p^{2k}-p^{2k-1}}) & \doteq  
\lambda_2 \cdot \mu^{-p^{2k-1}} \cdot (t\mu)^{r(2k)}\\
\end{split}
\end{equation*}
with $v_p(j)=2k-3$, leaving
\begin{equation*}
\begin{split}
\mu^{-1} E^{2r(2k)+1}(C_{p^n})=\ &
P(\mu^{\pm p^{2k}}) \otimes E(u_n, \lambda_1,a_1) \otimes P_{p-1}(b_1)
\otimes P(t\mu)\\
& \oplus E(u_n, \lambda_1) \otimes P_{p-2}(b_1) \otimes P(t\mu)
\otimes \F_p\{a_1\mu^j,b_1\mu^j\,|\,v_p(j)=2k-1\}\\
& \oplus \bigoplus_{1\leq m\leq k} T_m(C_{p^n}),
\end{split}
\end{equation*}
where for $m\geq2$ we have
\begin{equation*}
\begin{split}
 T_m(C_{p^n})=\ &
E(u_n,\lambda_1) \otimes P_{r(2m)}(t\mu)\otimes
\F_p\{\lambda_2\mu^j\,|\,v_p(j)=2m-1\}\\
&\oplus E(u_n,a_1) \otimes P_{p-1}(b_1) \otimes P_{r(2m-1)}(t\mu)\otimes
\F_p\{\lambda_1\mu^j\,|\,v_p(j)=2m-2\}\\
&\oplus E(u_n,\lambda_1) \otimes P_{p-2}(b_1) \otimes P_{r(2m-2)+1}(t\mu)\otimes
\F_p\{a_1\mu^j\,|\,v_p(j)=2m-3\}\,.
\end{split}
\end{equation*}
For $n\geq1$, there is a last non-zero differential
\begin{equation*}
d^{2r(2n)+1}(u_n\cdot \mu^{p^{2n}})  \doteq  (t\mu)^{r(2n)+1}
\end{equation*}
after which the spectral sequence collapses, leaving
\begin{equation*}
\begin{split}
\mu^{-1} E^{\infty}(C_{p^n}) =\ & 
P(\mu^{\pm p^{2n}}) \otimes E(\lambda_1,a_1) \otimes P_{p-1}(b_1)
\otimes P_{r(2n)+1}(t\mu)\\
& \oplus E(\lambda_1) \otimes P_{p-2}(b_1) \otimes P_{r(2n)+1}(t\mu)
\otimes
\F_p\{a_1\mu^j,b_1\mu^j\,|\,v_p(j)=2n-1\}\\
& \oplus \bigoplus_{1\leq m\leq n} T_m(C_{p^n})\,.
\end{split}
\end{equation*}
\end{lemma}
\begin{proof}
We prove these two lemmas by induction on $n$, showing that
Lemma~\ref{kku-lem:cpn-tate} for $C_{p^n}$ implies 
Lemma~\ref{kku-lem:cpn-fixed} for $C_{p^n}$, which in turn implies
Lemma~\ref{kku-lem:cpn-tate} for $C_{p^{n+1}}$.
The induction starts with Lemma~\ref{kku-lem:cpn-tate} for $C_p$, 
which is the content of Lemma~\ref{kku-lem:cp-tate}.
Let us therefore assume given $n\geq1$ such that 
Lemma~\ref{kku-lem:cpn-tate} holds for
$C_{p^n}$. 
The homotopy restriction map
\begin{equation*}
R^h:\thh(ku_p)^{hC_{p^n}}\to\thh(ku_p)^{tC_{p^n}} 
\end{equation*} 
induces a morphism of
spectral sequences $(R^h)^*:E^*(C_{p^n})\to
\hat E^*(C_{p^n})$, which at the $E^2$-terms corresponds to 
inverting the class $t\in E^2_{-2,0}(C_{p^n})$,
\begin{equation}\label{kku-eq:rh2} 
(R^h)^2 : E^2(C_{p^n}) \subset E^2(C_{p^n})[t^{-1}]\cong 
\hat E^2(C_{p^n})\,,
\end{equation}
and can be pictured as the inclusion of the second
quadrant into the upper-half plane. 
As we will see below, although $(R^h)^r$ is not injective
for $r\geq3$, it detects all the non-trivial differentials of
$E^r(C_{p^n})$. Taking into account the multiplicative
structure and the fact that $\lambda_1, \lambda_2, b_1$ and $t\mu$ are
infinite cycles, we claim that these differential are given by 
\begin{equation*}
\begin{split}
d^2(b_i)&=(1-i)a_i t\\
d^{2p}(t)&\doteq \lambda_1\cdot t^{1+p}\\
d^{2p^2}(t^p)&\doteq \lambda_2\cdot t^{p+p^2}
\end{split}
\end{equation*} 
with $0\leq i\leq p-1$, 
\begin{equation*}
\begin{split}
d^{2r(2k)+2}(b_1t^j) & \doteq 
a_1t^j \cdot t^{p^{2k}} \cdot (t\mu)^{r(2k-2)+1}\\
d^{2r(2k+1)}(t^{p^{2k}}) & \doteq  
\lambda_1 \cdot t^{p^{2k}+p^{2k+1}} \cdot (t\mu)^{r(2k-1)}\\
d^{2r(2k+2)}(t^{p^{2k+1}}) & \doteq  
\lambda_2 \cdot t^{p^{2k+1}+p^{2k+2}} \cdot (t\mu)^{r(2k)}\\
\end{split}
\end{equation*}
if $n\geq2$, $1\leq k\leq n-1$ and $v_p(j)=2k-1$ with $i\geq0$, and finally 
\begin{equation*}
d^{2r(2n)+1}(u_n)  \doteq  (t\mu)^{r(2n-2)+1}\cdot
t^{p^{2n}}\,.
\end{equation*}
To prove this claim, we assume that
some $r\geq2$ is given, and that $E^r(C_{p^n})$ has been
computed using the differentials $d^{r'}$ above with $r'<r$. 
The class $t\mu$ is an infinite cycle, and 
$E^r(C_{p^n})$ is a $P(t\mu)$-module. Our choice of generators induces
a decomposition $E^r(C_{p^n})\cong F^r(C_{p^n})\oplus T^r(C_{p^n})$,
where  $F^r(C_{p^n})$ is a free $P(t\mu)$-module and
$T^r(C_{p^n})$ is a $t\mu$-torsion module. 
By inspection, the non-zero elements of $T^r(C_{p^n})$ are
concentrated in filtration degrees $s$ with $-r<s\leq0$, so they
cannot be boundaries. They cannot support non-zero differentials either since
a $t\mu$-torsion class cannot map to a non-torsion class. Thus
the differential $d^r$ maps $F^r(C_{p^n})$ to itself and
$T^r(C_{p^n})$ to zero. The 
morphism $(R^h)^r$ maps $F^r(C_{p^n})$
injectively into $\hat E^r(C_{p^n})$, and it therefore detects
the non-zero differentials of $E^r(C_{p^n})$ as the 
non-zero differential of $\hat E^r(C_{p^n})$ which lie in
the second quadrant. These are precisely the differentials
given above. By induction on $r$, this determines all the non-trivial
differentials of $E^*(C_{p^n})$.
In the $\mu$-inverted homotopy fixed-point spectral sequence
$\mu^{-1}E^*(C_{p^n})$, these can be rewritten as the claimed
differentials. This proves Lemma~\ref{kku-lem:cpn-fixed} for $C_{p^n}$.
\par
We now turn to the proof of Lemma~\ref{kku-lem:cpn-tate} for $C_{p^{n+1}}$.
In the Tate spectral sequence $\hat E^*(C_{p^n})$ the first 
non-zero differential of odd length originating from a column of odd
$s$-filtration is $d^{2r(2n)+1}$. By~\cite{Rkku-AR02}*{Lemma~5.2} the spectral
sequences $\hat E^*(C_{p^n})$ and $\hat E^*(C_{p^{n+1}})$ are
abstractly isomorphic up to the $E^{2r(2n)+1}$-term included. The
Frobenius map 
\begin{equation*}
F:\thh(ku_p)^{tC_{p^{n+1}}}\to \thh(ku_p)^{tC_{p^{n}}}
\end{equation*}
induces a morphism of the corresponding Tate spectral sequences, which
on $E^r$-terms with $2\leq r\leq 2r(2n)+1$ maps the columns of even
$s$-filtration isomorphically. This detects all the claimed differentials of $\hat
E^r(C_{p^{n+1}})$ for $2\leq r\leq 2r(2n)$, and leaves
\begin{equation*}
\hat E^{2r(2n)+1}(C_{p^{n+1}})=\hat
F^{2r(2n)+1}(C_{p^{n+1}})\oplus\bigoplus_{m=1}^{n-1}\hat
T_m(C_{p^{n+1}})\,,
\end{equation*}
where $\hat F^{2r(2n)+1}(C_{p^{n+1}})$ is the $t\mu$-torsion free summand
\begin{equation*}
\begin{split}
\hat F&^{2r(2n)+1}(C_{p^{n+1}}) =  
P(t^{\pm p^{2n}}) \otimes E(u_{n+1},\lambda_1) \otimes P(t\mu) \\
& \otimes \Big( P_{p-1}(b_1) \otimes E(a_1) 
\oplus  P_{p-2}(b_1) \otimes \F_p\{a_1t^{-ip^{2n-1}},b_1t^{-ip^{2n-1}}
\,|\,0<i<p \} \Big)\,.
\end{split}
\end{equation*}
The non-zero $t\mu$-torsion elements of $\hat
E^{2r(2n)+1}(C_{p^{n+1}})$ are concentrated in internal
degrees $t$ with $0\leq t< 2r(2n)$.
In particular these elements cannot be boundaries, and they cannot
map to non-$t\mu$-torsion elements. As in the case of the homotopy
fixed-point spectral sequence above, we deduce that for $r\geq
2r(2n)+1$ the differential
$d^r$ can only affect the summand $\hat
F^{2r(2n)+1}(C_{p^{n+1}})$.
By Lemma~\ref{kku-lem:ss-desc} the summand of $\delta$-weight $0$ of
$\hat E^*(C_{p^{n+1}})$ is equal to the image of the
injective morphism of spectral sequences
\begin{equation*}
\hat E^*(C_{p^{n+1}},\ell_p)\to  \hat E^*(C_{p^{n+1}},ku_p)= \hat
E^*(C_{p^{n+1}})
\end{equation*}
induced by the map $\ell_p\to ku_p$.
Therefore, by~\cite{Rkku-AR02}*{Theorem~6.1}, the differentials affecting the
summand of $\delta$-weight $0$ of $\hat F^{2r(2n)+1}(C_{p^{n+1}})$ at
a later stage are given by
\begin{equation}\label{kku-eq:diff-tate}
\begin{split}
d^{2r(2n+1)}(t^{p^{2n}-p^{2n+1}}) & \doteq  
\lambda_1 \cdot t^{p^{2n}} \cdot (t\mu)^{r(2n-1)}\\
d^{2r(2n+2)}(t^{p^{2n+1}-p^{2n+2}}) & \doteq  
\lambda_2 \cdot t^{p^{2n+1}} \cdot (t\mu)^{r(2n)} \\
d^{2r(2n+2)+1}(u_{n+1}\cdot t^{-p^{2n+2}}) & \doteq
(t\mu)^{r(2n)+1}\,,
\end{split}
\end{equation} 
together with the multiplicative structure and the fact that $t\mu$
is an infinite cycle. It remains to prove that from the 
$E^{2r(2n)+1}$-term on, the
only non-zero differentials supported by homogeneous algebra generators of
$\delta$-weight $1$ are given by
\begin{equation}\label{kku-eq:dif-tat2}
d^{2r(2n)+2}(b_1t^j) \doteq 
a_1t^j \cdot t^{p^{2n}} \cdot (t\mu)^{r(2n-2)+1}
\end{equation}
for $v_p(j)=2n-1$. 
First, notice that for tri-degree reasons $d^{2r(2n)+1}=0$, so that
$\hat F^{2r(2n)+2}(C_{p^{n+1}}) = \hat
F^{2r(2n)+1}(C_{p^{n+1}})$.
To detect the differential~\eqref{kku-eq:dif-tat2} we make use
of the $(2p-2)$-coconnected map 
\begin{equation*}
  (\hat\Gamma_{n+1})_*:V(1)_*\thh(ku_p)^{C_{p^n}}\to
V(1)_*\thh(ku_p)^{tC_{p^{n+1}}}\,, 
\end{equation*}
and argue as
in~\cite{Rkku-AR02}*{proof of~6.1}. 
There is a commutative diagram
\begin{equation*}
\xymatrix{
\thh(ku_p)^{hC_{p^n}} \ar[d]^{F^n} &
\thh(ku_p)^{C_{p^n}}  \ar[d]^{F^n} \ar[l]_-{\Gamma_n} \ar[r]^-{\hat\Gamma_{n+1}}&
\thh(ku_p)^{tC_{p^{n+1}}} \ar[d]^{F^n} 
\\
\thh(ku_p)  &
\thh(ku_p)  \ar[l]_-{\Gamma_0}^= \ar[r]^-{\hat\Gamma_1}&
\thh(ku_p)^{tC_{p}}
}
\end{equation*}
where the vertical arrows are the $n$-fold Frobenius maps. 
The left-hand Frobenius is given in $V(1)$-homotopy on the associated graded
by the edge homomorphism 
\begin{equation*}
E^\infty_{*,*}(C_{p^n}) \to E^\infty_{0,*}(C_{p^n})\subset
E^2_{0,*}(C_{p^n})=V(1)_*\thh(ku_p)\,,
\end{equation*}
which is known by induction hypothesis. 
For each $0<\ell <p$ there is a direct summand 
\begin{equation*}
P_{r(2n-2)+1}(t\mu)\{a_1\mu^{\ell p^{2n-3}}\}\subset
E^\infty_{*,*}(C_{p^n})\,,
\end{equation*} 
and $a_1\mu^{\ell p^{2n-3}}$ maps by $F^n_*$ to the class with same name in
$V(1)_*\thh(ku_p)$. 
Since $(\Gamma_n)_*$ is $(2p-2)$-coconnected, there is
a class $x_\ell \in V(1)_*\thh(ku_p)^{C_{p^n}}$
with $F^n_*(x_\ell )=a_1\mu^{\ell p^{2n-3}}$ in $V(1)_*\thh(ku_p)$.
In $E^\infty(C_{p^n})$ we have no non-zero class of
same total degree, same $\delta$-weight and lower
$s$-filtration than 
\begin{equation*}
(t\mu)^{r(2n-2)+1}\cdot
a_1\mu^{\ell p^{2n-3}}\,,
\end{equation*}
which forces $v_2^{r(2n-2)+1}x_\ell =0$ in
$V(1)_*\thh(ku_p)^{C_{p^n}}$.
By Proposition~\ref{kku-prop:gamma1}, the class $(\hat\Gamma_1 F^n)_*(x_\ell )$ 
is represented by $a_1t^{-\ell p^{2n-1}}\in \hat E^\infty(C_p)$, 
and therefore $(\hat\Gamma_{n+1})_*(x_\ell )$ must be detected in 
$s$-filtration $2\ell p^{2n-1}$ or higher. 
The only suitable class in $\hat E^{2r(2n)+2}(C_{p^{n+1}})$  
is $a_1t^{-\ell p^{2n-1}}$, which therefore
is a permanent cycle
representing $(\hat\Gamma_{n+1})_*(x_\ell )$. 
Notice for later use that the same argument shows that
\begin{equation*}
a_1\in \hat E^{2r(2n)+2}_{0,2p+3}(C_{p^{n+1}})
\end{equation*} 
is a permanent cycle.  The map
$(\hat\Gamma_{n+1})_*$ is an isomorphism 
in degrees larger than $2p-2$, and the relation
$v_2^{r(2n-2)+1}(\hat\Gamma_{n+1})_*(x_\ell )=0$ implies that the
infinite cycle $(t\mu)^{r(2n-2)+1}\cdot a_1t^{-\ell
p^{2n-1}}$, of total degree $2p^{2n}+2\ell p^{2n-1}+2p+1$
and of $\delta$-weight~$1$, is a boundary. 
On the other hand, 
the component of $\hat E^{2r(2n)+1}(C_{p^{n+1}})$ of total
degree $2p^{2n}+2\ell p^{2n-1}+2p+2$, of $\delta$-weight~$1$ and
of $s$-filtration degree exceeding by at least $2r(2n)+2$
the $s$-filtration degree of $(t\mu)^{r(2n-2)+1}\cdot
a_1t^{-\ell p^{2n-1}}$ reduces to 
\begin{equation*}
\F_p\{b_1t^{-\ell p^{2n-1}}\cdot t^{-p^{2n}}\}\,. 
\end{equation*}
This proves the existence of a non-zero differential
\begin{equation*}
d^{2r(2n)+2}(b_1t^{-\ell p^{2n-1}}\cdot
t^{-p^{2n}})\doteq (t\mu)^{r(2n-2)+1}\cdot
a_1t^{-\ell p^{2n-1}}
\end{equation*}
for $0<\ell <p$. Since $t^{p^{2n}}$ is a unit and a
cycle we obtain the claimed differentials~\eqref{kku-eq:dif-tat2}.
This leaves
\begin{equation*}
\begin{split}
\hat E^{2r(2n)+3}(C_{p^{n+1}})&= \hat
F^{2r(2n)+3}(C_{p^{n+1}})\\
&\,\oplus E(u_{n+1},\lambda_1)\otimes P_{p-2}(b_1)\otimes
P_{r(2n-2)+1}(t\mu)\otimes \F_p\{a_1t^j\,|\,v_p(j)=2n-1 \}\\
&\,\oplus\bigoplus_{m=1}^{n-1}\hat
T_m(C_{p^{n+1}})\,,
\end{split}
\end{equation*} 
with a $t\mu$-torsion free summand 
\begin{equation*}
F^{2r(2n)+3}(C_{p^{n+1}})=P(t^{\pm p^{2n}})\otimes
E(u_{n+1}, \lambda_1, a_1)\otimes P_{p-1}(b_1)\otimes
P(t\mu)\,.
\end{equation*}
Again, further differentials can only affect the summand
$F^{2r(2n)+3}(C_{p^{n+1}})$. Since $b_1$ and $a_1$ are
infinite cycles, the next non-zero differentials are $d^{2r(2n+1)}$ and
$d^{2r(2n+2)}$, as given in~\eqref{kku-eq:diff-tate}, leaving
\begin{equation*}
\hat E^{2r(2n+2)+1}(C_{p^{n+1}})=\hat
F^{2r(2n+2)+1}(C_{p^{n+1}})\oplus\bigoplus_{m=1}^{n}\hat
T_m(C_{p^{n+1}})\,,
\end{equation*}
with 
\begin{equation*}
\begin{split}
\hat F&^{2r(2n+2)+1}(C_{p^{n+1}}) =  
P(t^{\pm p^{2n+2}}) \otimes E(u_{n+1},\lambda_1) \otimes P(t\mu) \\
& \otimes \Big( P_{p-1}(b_1) \otimes E(a_1) 
\oplus  P_{p-2}(b_1) \otimes \F_p\{a_1t^{-ip^{2n+1}},b_1t^{-ip^{2n+1}}
\,|\,0<i<p \} \Big)\,.
\end{split}
\end{equation*}
Notice that for tri-degree reasons,
the classes $a_1t^{-ip^{2n+1}}$ and
$b_1t^{-ip^{2n+1}}$ are cycles at the $E^{2r(2n+2)+1}$-stage.
The third differential of~\eqref{kku-eq:diff-tate} remains,
after which the spectral sequence collapses for bidegree
reasons, leaving
\begin{equation*}
\begin{split}
\hat E^{\infty}(C_{p^{n+1}}) =\ &
P(t^{\pm p^{2(n+1)}}) \otimes E(\lambda_1,a_1) \otimes P_{p-1}(b_1)
\otimes P_{r(2n)+1}(t\mu)\\
&\oplus E(\lambda_1) \otimes P_{p-2}(b_1) \otimes P_{r(2n)+1}(t\mu)
\otimes
\F_p\{a_1t^j,b_1t^j\,|\,v_p(j)=2n+1\}\\
& \oplus \bigoplus_{1\leq m\leq n} \hat T_m(C_{p^{n+1}})\,,
\end{split}
\end{equation*}
as claimed.
This completes the induction step and the proof of 
Lemmas~\ref{kku-lem:cpn-tate} and~\ref{kku-lem:cpn-fixed}. 
\end{proof}
\par
Taking the limit over the Frobenius maps we obtain the
following two lemmas.
\begin{lemma}\label{kku-lem:tate-S1}
The associated graded $\hat E^{\infty}(S^1)$ of
$V(1)_*\thh(ku_p)^{tS^1}$ is given by
\begin{equation*}
\hat E^{\infty}(S^1)=
E(\lambda_1,a_1) \otimes P_{p-1}(b_1) \otimes P(t\mu)
\oplus \bigoplus_{m\geq 1} \hat T_m(S^1)\,,
\end{equation*}
where
\begin{equation*} 
\begin{split}
\hat T_m(S^1) =\ &
E(\lambda_1) \otimes P_{r(2m)}(t\mu)\otimes
\F_p\{\lambda_2t^j\,|\,v_p(j)=2m+1\}\\
&\oplus E(a_1) \otimes P_{p-1}(b_1) \otimes P_{r(2m-1)}(t\mu)\otimes
\F_p\{\lambda_1t^j\,|\,v_p(j)=2m\}\\
&\oplus E(\lambda_1) \otimes P_{p-2}(b_1) \otimes P_{r(2m-2)+1}(t\mu)\otimes
\F_p\{a_1t^j\,|\,v_p(j)=2m-1\}\,.
\end{split}
\end{equation*}
\end{lemma}
\begin{lemma}\label{kku-lem:fixed-S1}
The associated graded $E^{\infty}(S^1)$ of
$V(1)_*\thh(ku_p)^{hS^1}$ is mapped by a $(2p-2)$-coconnected
homomorphism to
\begin{equation*}
\mu^{-1} E^{\infty}(S^1)=
E(\lambda_1,a_1) \otimes P_{p-1}(b_1) \otimes P(t\mu)
\oplus \bigoplus_{m\geq 1} T_m(S^1)\,,
\end{equation*}
where
\begin{equation*}
\begin{split}
T_1(S^1) =\ &
E(\lambda_1) \otimes P_{p^2}(t\mu)\otimes
\F_p\{\lambda_2\mu^j\,|\,v_p(j)=1\}\\
&\oplus E(a_1) \otimes P_{p-1}(b_1) \otimes P_{p}(t\mu)\otimes
\F_p\{\lambda_1\mu^j\,|\,v_p(j)=0\}\\
&\oplus E(\lambda_1) \otimes P_{p-2}(b_1) \otimes P(\mu^{\pm1})\otimes
\F_p\{a_i\,|\,0\leq i\leq p-1,\ i\neq1\}
\end{split}
\end{equation*}
and, for $m\geq2$,
\begin{equation*}
\begin{split}
 T_m(S^1)=\ &
E(\lambda_1) \otimes P_{r(2m)}(t\mu)\otimes
\F_p\{\lambda_2\mu^j\,|\,v_p(j)=2m-1\}\\
&\oplus E(a_1) \otimes P_{p-1}(b_1) \otimes P_{r(2m-1)}(t\mu)\otimes
\F_p\{\lambda_1\mu^j\,|\,v_p(j)=2m-2\}\\
&\oplus E(\lambda_1) \otimes P_{p-2}(b_1) \otimes P_{r(2m-2)+1}(t\mu)\otimes
\F_p\{a_1\mu^j\,|\,v_p(j)=2m-3\}\,.
\end{split}
\end{equation*}
\end{lemma}
\section{Topological cyclic homology}
\label{kku-sec:tc}
We now evaluate the restriction map $R:TF(ku_p)\to TF(ku_p)$ in
$V(1)$-homo\-to\-py. 
Consider the homotopy commutative diagram
\begin{equation*}
\xymatrix{
TF(ku_p)  \ar[r]^-R  \ar[d]^{\Gamma} &
TF(ku_p)  \ar[d]^{\hat\Gamma} \ar[r]^-{\Gamma} &
\thh(ku_p)^{hS^1} \ar[d]^{({\hat\Gamma}_1)^{hS^1}} \\
\thh(ku_p)^{hS^1}  \ar[r]^-{R^h} &
\thh(ku_p)^{tS^1}  \ar[r]^-G &
(\thh(ku_p)^{tC_p})^{hS^1}                       
}
\end{equation*}
displayed in~\cite{Rkku-AR02}*{page~27}, and with $G$ a 
$V(1)$-equivalence.
By the argument in~\cite{Rkku-AR02}*{Lemma 7.5}, we know
that on $V(1)_*TF(ku_p)$ the profinite topology coincides with the topology induced
by the spectral sequence filtration of
$V(1)_*\thh(ku_p)^{hS^1}$ via $\Gamma_*$, and that the restriction map
\begin{equation*}
  R_*:V(1)_*TF(ku_p)\to V(1)_*TF(ku_p)
\end{equation*}
is continuous in degrees larger than $2p-2$.
In this range of degrees, we will identify
$V(1)_*TF(ku_p)$ with $V(1)_*\thh(ku_p)^{hS^1}$ via the
homeomorphism $\Gamma_*$. Under this
identification $R_*$ corresponds to
$(\Gamma_*\hat\Gamma_*^{-1}) R^h_*$, and we 
first describe $R^h_*$ and $\Gamma_*\hat\Gamma_*^{-1}$
separately.
\begin{lemma}\label{kku-lem:rhinfty}
In total degrees larger than $2p-2$, the morphism
\begin{equation*}
(R^h)^\infty:E^\infty(S^1)\to \hat E^\infty(S^1)
\end{equation*}
has the following properties.
\begin{itemize}
\item[\textup{(a)}] It maps $E(\lambda_1,a_1) \otimes P_{p-1}(b) \otimes
P(t\mu)$ isomorphically to the  summand with same name\,;
\item[\textup{(b)}] It maps $E(\lambda_1)\otimes P_{r(k)}(t\mu)\otimes
\F_p\{\lambda_2\mu^{-dp^{k-1}}\}$ onto 
\begin{equation*}
E(\lambda_1)\otimes P_{r(k-2)}(t\mu)\otimes
\F_p\{\lambda_2 t^{dp^{k-1}}\} 
\end{equation*}
and 
$E(\lambda_1)\otimes P_{p-2}(b)\otimes P_{r(k)+1}(t\mu)\otimes
\F_p\{a_1\mu^{-dp^{k-1}}\}$ onto 
\begin{equation*}
E(\lambda_1)\otimes P_{p-2}(b)\otimes P_{r(k-2)+1}(t\mu)\otimes
\F_p\{a_1 t^{dp^{k-1}}\} 
\end{equation*}
for $k\geq2$ even and $0<d<p$\,; 
\item[\textup{(c)}] It maps
$E(a_1)\otimes P_{p-1}(b)\otimes P_{r(k)}(t\mu)\otimes
\F_p\{\lambda_1\mu^{-dp^{k-1}}\}$ onto 
\begin{equation*}
E(a_1)\otimes P_{p-1}(b)\otimes P_{r(k-2)}(t\mu)\otimes
\F_p\{\lambda_1 t^{dp^{k-1}}\}
\end{equation*}
for $k\geq3$ odd and $0<d<p$\,; 
\item[\textup{(d)}] It maps the remaining summands to zero.
\end{itemize}
\end{lemma}
\begin{proof}
This follows from the description of $(R^h)^2$,
see~\eqref{kku-eq:rh2}.
\end{proof}
\begin{lemma}\label{kku-lem:r-gamma}
In degrees larger then $2p-2$, the homomorphism
$\Gamma_*\hat\Gamma_*^{-1}$ maps
\begin{itemize}
\item[\textup{(a)}] the classes in $V(1)_*\thh(ku_p)^{tS^1}$ represented 
in $\hat E^\infty(S^1)$ by 
\begin{equation*}
\lambda_1^{\epsilon_1} a_1^{\epsilon_2} b^k (t\mu)^m t^i
\end{equation*}
for $v_p(i)\neq 1$, $\epsilon_1$ and $\epsilon_2\in\{0,1\}$,
$0\leq k\leq p-2$ and $m\geq0$, to
classes in \\
$V(1)_*\thh(ku_p)^{hS^1}$ represented in
$E^\infty(S^1)$  by
\begin{equation*}
\lambda_1^{\epsilon_1} a_1^{\epsilon_2} b^k (t\mu)^m
\mu^j
\end{equation*} 
with $i+p^2j=0$, 
up to multiplication with a unit in $\F_p$\,; 
\item[\textup{(b)}] the classes in $V(1)_*\thh(ku_p)^{tS^1}$ represented in
$\hat E^\infty(S^1)$ by
\begin{equation*}
\lambda_1^{\epsilon_1}b^k a_1t^{i}
\end{equation*}
for $v_p(i)=1$, 
$\epsilon_1\in\{0,1\}$ and $0\leq k\leq p-3$,
to classes in $V(1)_*\thh(ku_p)^{hS^1}$ represented in
$E^\infty(S^1)$ by
\begin{equation*}
\lambda_1^{\epsilon_1} b^k \mu^l a_j
\end{equation*} 
with 
$i=(1-j)p-lp^2$ 
for 
$0\leq j\leq p-1$ such that $j\neq1$,
up to multiplication with a unit in $\F_p$.
\end{itemize}
\end{lemma}
\begin{proof}
The proof is similar to the proof
of~\cite{Rkku-AR02}*{Proposition~7.4}, and we omit it.
\end{proof}
\begin{definition}\label{kku-def:TF-classes}
We recall from~\cite{Rkku-AR02}*{Theorem~9.1} that there are
classes $\lambda_1t^{p-1}$, $\lambda_1$ and $\lambda_2$
in $V(1)_*K(\ell_p)\subset V(1)_*K(ku_p)$, of degree $1$, $2p-1$ and $2p^2-1$,
respectively. We denote by 
\begin{equation*}
\wti{\lambda_1t^{p-1}},\ \ \ti{\lambda}_1\
\textup{ and }\ \ti{\lambda}_2 
\end{equation*}
their image in $V(1)_*TF(ku_p)$ under ${\tr_F}_*$. The latter
classes are represented by
\begin{equation*}
\lambda_1t^{p-1}=(t\mu)^{p-1}\cdot\lambda_1\mu^{1-p},\
\lambda_1\ \textup{ and }\ \lambda_2 
\end{equation*}
in $E^\infty(S^1)$, respectively,
see~\cite{Rkku-AR02}*{Theorem~8.4}. 
We further denote by $b$ and $v_2$ the image in
$V(1)_*TF(ku_p)$ under ${\tr_F}_*$ of the classes with same
name in $V(1)_*K(ku_p)$. These classes are represented
by $b_1$ and $t\mu$ in $E^\infty(S^1)$, respectively, see
Lemma~\ref{kku-lem:b->b1} and~\cite{Rkku-AR02}*{Proposition~4.8}.
\end{definition}
\begin{lemma}\label{kku-lem:tilde-a1}
There exists a unique class $\ti{a}_1\in V(1)_{2p+3}TF(ku_p)$
with the following two properties\,:
\begin{itemize}
  \item[\textup{(a)}] $\ti{a}_1$ has $\delta$-weight $1$ and 
    $b^{p-2} \ti{a}_1=\ti{\lambda}_2$,
  \item[\textup{(b)}] $R_*(\ti{a}_1)=\ti{a}_1$.
\end{itemize}
Moreover, this class $\ti{a}_1$ is represented by $a_1$ in
    $E^\infty(S^1)$. 
\end{lemma}
\begin{proof}
For $i=0$ or $1$, let us denote by $T_*^{(i)}$ and
$\ker(R-1)_*^{(i)}$ the
summand of $\delta$-weight $i$ of $V(1)_*TF(ku_p)$ and
$\ker(R-1)_*\subset V(1)_*TF(ku_p)$, 
respectively.
We make the following claims\,:
\begin{itemize}
  \item[(1)] The homomorphism given by multiplication with
    $b^{p-2}$ on $T_{2p+3}^{(1)}$ fits in a short exact sequence
    \begin{equation*}
      0\to \F_p\{z\}\to
      T^{(1)}_{2p+3}\stackrel{b^{p-2}}{\lr}
      T^{(0)}_{2p^2-1}\to 0\,,
    \end{equation*} 
    where the class $z$ is represented by
    $b_1\cdot(t\mu)^{p-1}\cdot \lambda_1\mu^{1-p}$ in
    $E^\infty(S^1)$\,;
  \item[(2)] The class $z$ does not belong to $\ker(R-1)_*$.
\end{itemize}
Using these claims, it is easy to deduce that multiplication
with $b^{p-2}$ restricts to an isomorphism
\begin{equation*}
 \ker(R-1)_{2p+3}^{(1)}\stackrel{\cong}{\lr}\ker(R-1)_{2p^2-1}^{(0)}\,.
\end{equation*}  
We have ${\ti\lambda}_2\in\ker(R-1)_{2p^2-1}^{(0)}$ since
${\ti\lambda}_2$ has $\delta$-weight $0$ and is in the image of ${\tr_F}_*$.
Therefore, there is a unique pre-image ${\ti a}_1\in
\ker(R-1)_{2p+3}^{(1)}$ of
${\ti\lambda}_2\in\ker(R-1)_{2p^2-1}^{(0)}$, or, in other
words, there is a unique class ${\ti a}_1\in
V(1)_{2p+3}TF(ku_p)$ with properties (a) and (b).  
Moreover, ${\ti\lambda}_2$ is
represented in $E^\infty(S^1)$ in filtration zero by $\lambda_2=b_1^{p-2}a_1$,
and we deduce that ${\ti a}_1$ must be represented in
filtration zero by $a_1$. Thus this lemma follows 
from claims (1) and (2), which we now prove. 
\par
First, notice that 
the group $T_*^{(i)}$ inherits via $\Gamma_*$ the
spectral sequence filtration of $V(1)_*\thh(ku_p)^{hS^1}$.
Denoting by $E^\infty(S^1)^{(i)}_*$ its associated
graded, we know from Lem\-ma~\ref{kku-lem:fixed-S1} that
\begin{equation*}
  \begin{split}
    E^\infty(S^1)^{(1)}_{2p+3} =\ &  
    \F_p\{a_1, b_1\cdot x_n\,|\,n\geq 0\} \textup{ and }\\
    E^\infty(S^1)^{(0)}_{2p^2-1} =\ & 
    \F_p\{\lambda_2, t\mu\cdot x_n\,|\,n\geq 1\}\,,\\
  \end{split}
\end{equation*}
where
$x_n=(t\mu)^{r(2n+1)-r(2n)-1}\cdot\lambda_1\mu^{(1-p)p^{2n}}$.
\par
Next, the relation
$b'^p+v_2b'=0$ in $V(1)_*K(\Z,3)$, established
in Proposition~\ref{kku-prop:prebott}, maps under ${\tr_F}_*\phi_*$
to the relation $b^p+v_2b=0$ in $T^{(1)}_*$.
The class $v_2b$ in $T^{(1)}_*$ is represented by the non-zero
class $t\mu\cdot b_1$ in $E^\infty(S^1)$ in filtration $-2$, 
and we deduce that $b^{p-1}\in T^{(0)}_*$
must be represented by $-t\mu$ in $E^\infty(S^1)$. It
follows that if a class $x\in T^{(1)}_{2p+3}$ is represented
by $b_1\cdot x_n$, then $b^{p-2}x$ is represented by
$-t\mu\cdot x_n$ in $2$ filtration degrees lower. Using a coarser
filtration that ignores this shift, 
and considering our formulas for $E^\infty(S^1)^{(1)}_{2p+3}$
and $E^\infty(S^1)^{(0)}_{2p^2-1}$ given above, 
we deduce claim (1) from the corresponding claim 
for the associated graded, with $z$ represented by $b_1\cdot
x_0$.
\par
\smallskip
To prove claim (2), we notice that if a class $y\in
T_{2p+3}^{(1)}$ is represented by $b_1\cdot x_n$ with $n\geq1$, then 
$R_*(y)$ will be represented by $b_1\cdot x_{n-1}$ in higher
filtration, up to some non-zero scalar multiple\,: this
follows directly from Lemmas~\ref{kku-lem:rhinfty}
and~\ref{kku-lem:r-gamma}. In particular, $R_*(y)\neq y$.
This implies the following claim\,:
\begin{itemize}
  \item[(3)] The group $\ker(R-1)_{2p+3}^{(1)}$ contains at most one
    class represented by $b_1\cdot x_0$.
\end{itemize}
Now consider the class $\ti x_0=\wti{\lambda_1t^{p-1}}\in
T^{(0)}_1$ given in Definition~\ref{kku-def:TF-classes}. 
By definition, this class lies in $\ker(R-1)^{(0)}_1$ and is 
represented by $x_0$. We also claim that
\begin{itemize}
  \item[(4)] The class $b\ti x_0\in \ker(R-1)_{2p+3}^{(1)}$ is not
    annihilated by $b^{p-2}$.
\end{itemize}
Since $b\ti x_0$ is represented by $b_1\cdot x_0$, claim~(2)
follows from claims~(3) and~(4). 
\par
\smallskip
Finally, to prove claim~(4), we recall
from~\cite{Rkku-AR02}*{Theorem~8.2}
that the class $v_2\ti x_0\in \ker(R-1)_{2p^2-1}^{(0)}$ is non-zero, and
must be represented, in filtration degree lower
then $-2p+2$, by a class in
\begin{equation*}
  \F_p\{t\mu\cdot x_n\,|\, n\geq 1 \}\,.
\end{equation*}
None of these classes is annihilated by $b_1$.
Therefore $bv_2\ti x_0=-b^p\ti x_0$ is non-zero, and we deduce that
$b\ti x_0\in\ker(R-1)_{2p+3}^{(1)}$ is not annihilated by $b^{p-2}$. 
\end{proof}
\smallskip
\begin{remark}
The lemma above implies that the class
$a_1\in V(1)_*\thh(ku_p)$ has a lift $a_1\in
V(1)_*K(ku_p)$ under the trace, with $b^{p-2}a_1=\lambda_2$, 
see Theorem~\ref{kku-thm:kku}. It would be nice 
to have a more direct construction of such a lift. In 
fact, we conjecture that $a_1\in V(1)_*K(ku_p)$ decomposes as $bd$, where
$d\in V(1)_1K(KU_p)$ is a unit class, when mapped into
$V(1)_*K(KU_p)$, see the discussion preceding
Theorem~\ref{kku-thm:periodic} below. 
\end{remark}
\begin{definition}\label{kku-def:some-k-classes}
We consider the following subgroups of $E^\infty(S^1)$\,: 
\begin{equation*}
\begin{split}
A  =\ & E(\lambda_1,a_1)\otimes P_{p-1}(b_1)\otimes
P(t\mu)\,,\\
B_0 =\ & 
E(\lambda_1)\otimes P_{p-2}(b_1)\otimes
\F_p\{\mu^{-1}a_i,
a_0\,|\, 2\leq i\leq p-1\}\,,\\
B_k =\ & 
\Big(E(\lambda_1)\otimes P_{p-2}(b_1)\otimes
\bigoplus_{0<d<p}\big(P_{r(k)-dp^{k-1}+1}(t\mu)\otimes
\F_p\{a_1t^{dp^{k-1}}\}\big)\Big)\\
&\ \oplus 
\Big(E(\lambda_1)\otimes 
\bigoplus_{0<d<p}P_{r(k)-dp^{k-1}}(t\mu)\otimes
\F_p\{\lambda_2t^{dp^{k-1}}\}\Big)
\ \hbox{\ for $k\geq2$ even,}\\
B_k =\ & 
E(a_1)\otimes P_{p-1}(b_1)\otimes\bigoplus_{0<d<p}
\big(P_{r(k)-dp^{k-1}}(t\mu)\otimes
\F_p\{\lambda_1t^{dp^{k-1}}\}\big)\ \hbox{\ for $k\geq1$
odd,}\\
\end{split}
\end{equation*}
and we let $C$ be the span of the remaining monomials in
$E^\infty(S^1)$. We then have a direct sum decomposition
$E^\infty(S^1)=A\oplus B\oplus C$, with
$B=\bigoplus_{k\geq0}B_k$. 
\end{definition}
\begin{lemma}\label{kku-lem:restr}
In dimensions larger than $2p-2$ there are closed subgroups
$\ti A$, $\ti B_k$ and $\ti C$ in $V(1)_*TF(ku_p)$,
represented by
$A$, $B_k$ and $C$ in $E^\infty(S^1)$ respectively, such that
\begin{itemize}
\item[\textup{(a)}] $R_*$ restricts to the identity on $\ti A$,
\item[\textup{(b)}] $R_*$ maps $\ti B_{k+2}$ onto $\ti B_{k}$ for $k\geq0$, 
\item[\textup{(c)}] $R_*$ maps $\ti B_0$, $\ti B_1$ and $\ti C$ to zero.
\end{itemize}
In these degrees $V(1)_*TF(ku_p)\cong \ti A \oplus\ti B\oplus \ti C$, where
$\ti B=\prod_{k\geq0}\ti B_k$.
\end{lemma}
\begin{proof}
On the associated
graded $E^\infty(S^1)$, the homomorphism $(\Gamma_*\hat\Gamma^{-1}_*)R^h_*$ 
has been described in Lemmas~\ref{kku-lem:rhinfty}
and~\ref{kku-lem:r-gamma}, and maps $A$ isomorphically to itself,
$B_{k+2}$ onto $B_k$ for $k\geq0$, and $B_0$, $B_1$ and $C$
to zero. It remains to find closed lifts of these groups in 
$V(1)_*TF(ku_p)$ with desired properties. 
We take $\ti A$ to be the (closed) subalgebra of
$V(1)_*TF(ku_p)$ generated by $\ti{\lambda}_1$, $\ti{a}_1$, $b$ and
$v_2$. Then $\ti A$ lifts $A$, by definition of
its algebra generators and by the fact, proved above,
that $b^{p-1}$ is represented by $-t\mu$ in $E^\infty(S^1)$. Also,
$\ti\lambda_1$, $b$ and $v_2$ are fixed under $R_*$, since
they are in the image of ${\tr_F}_*$, and $\ti a_1$ is fixed
by definition. To construct $\ti B_k$ for $k\geq0$ and
$\ti C$, we follow the procedure given
in~\cite{Rkku-AR02}*{Theorem~7.7}. 
\end{proof}
\begin{definition}
We denote $b\in V(1)_{2p+2}TC(ku_p)$ the image of the higher
Bott element $b$, defined in~\ref{kku-def:bott}, under the cyclotomic trace map
\begin{equation*}
({\mathrm{trc}})_*:V(1)_* K(ku_p)\to V(1)_*TC(ku_p)\,.
\end{equation*}
\end{definition}
\begin{theorem}\label{kku-thm:tcku}
The class $b\in V(1)_{2p+2}TC(ku_p)$ satisfies the relation
\begin{equation*}
b^{p-1}=-v_2\,.
\end{equation*}
There is an isomorphism of $P(b)$-modules
\begin{equation*}
\begin{split}
V(1)_*TC(ku_p)\cong\ & P(b)\otimes E(\partial, \lambda_1, a_1)
\\
&\oplus  P(b)\otimes E(a_1)\otimes \F_p\{t^d\lambda_1\,|\,0<d<p \}\\
&\oplus  P(b)\otimes
E(\lambda_1)\otimes\F_p\{u^ia_0,t^{p^2-p}\lambda_2\,|\,0\leq
i<p-2 \}\,,
\end{split}
\end{equation*}
where the degree of the classes is  
$|\partial|=-1$, $|\lambda_1|=2p-1$, $|a_1|=2p+3$, $|u^ia_0|=2i+3$,
$|\lambda_2|=2p^2-1$ and $|t|=-2$.
\end{theorem}
\begin{proof}
Recall that $TC(ku_p)$ is defined as the homotopy fiber of the 
map 
\begin{equation*}
R-1:TF(ku_p)\to TF(ku_p)\,.
\end{equation*}
In $V(1)$-homotopy, it gives a short exact 
sequence of $P(v_2)$-modules
\begin{equation}\label{kku-eq:co-ker}
0\to\Sigma^{-1}\cok(R-1)_*\to V(1)_*TC(ku_p) \to
\ker(R-1)_*\to 0\,.
\end{equation}
We have isomorphisms of $P(v_2)$-modules
\begin{equation}\label{kku-eq:co-ker2}
\begin{split}
\Sigma^{-1}\cok(R-1)_*&\cong \Sigma^{-1}\ti A\\
\ker(R-1)_*&\cong \ti A \oplus \lim_{k\geq 0\
{\textup{even}}}\ti B_k 
\oplus \lim_{k\geq 1\ {\textup{odd}}}\ti B_k\,.
\end{split}
\end{equation}
Indeed, $R_*-1$ maps each factor of the decomposition
$V(1)_*TF(ku_p)\cong \ti A\oplus\ti B\oplus \ti C$ to itself.
It restricts to zero on $\ti A$ and to the identity on $\ti C$. 
We have a short exact sequence
\begin{equation*}
0\to \lim_{k\geq 0\ {\textup{even}}}\ti B_k\to 
\prod_{k\geq 0\ {\textup{even}}}\ti B_k \stackrel{R_*-1}{\lr} 
\prod_{k\geq 0\ {\textup{even}}}\ti B_k
\to \Rlim_{k\geq 0\ {\textup{even}}}\ti B_k
\to 0\,,
\end{equation*}
and similarly for the $\ti B_k$ with $k$ odd.
Here the limits are taken over the sequential system of maps
$R_*:\ti B_{k+2}\to \ti B_k$ for $k\geq 0$ even or $k\geq 1$ odd. 
Since these maps are surjective, the $\Rlim$-terms are
trivial. This proves our claims on $\Sigma^{-1}\cok(R-1)_*$
and $\ker(R-1)_*$ in~\eqref{kku-eq:co-ker2}.
\par
For $k\geq1$ odd, 
the group $\ti B_k$ is isomorphic as a $P(v_2)$-module to a sum
of $2(p-1)^2$ cyclic $P(v_2)$-modules
\begin{equation*}
\ti B^k\cong E(a_1)\otimes P_{r(k)}(v_2)\otimes
P_{p-1}(b)\otimes
\F_p\{\lambda_1t^{dp^{k-1}}\,|\,0<d<p\}\,.
\end{equation*}
The map $R_*$ respects this decomposition into cyclic
$P(v_2)$-modules. Since the height of these modules grows to
infinity with $k$, we deduce from the surjectivity of $R_*$
that $\lim_{k\geq 1\ {\textup{odd}}}\ti B_k$
is a sum of $2(p-1)^2$ free cyclic $P(v_2)$-modules,
given by an isomorphism 
\begin{equation*} 
\lim_{k\geq 1\ {\textup{odd}}}\ti B_k\cong E(a_1)\otimes
P(v_2)\otimes P_{p-1}(b)\otimes
\F_p\{\lambda_1t^{d}\,|\,0<d<p\}\,.
\end{equation*}
Similarly, for $k\geq 2$ even,
$\ti B_k$ is isomorphic to a sum
of $2(p-1)^2$ cyclic $P(v_2)$-modules of height growing with
$k$, and passing to the limit we have an isomorphism of
$P(v_2)$-modules
\begin{equation*} 
\lim_{k\geq 0 \ {\textup{even}}}\ti B_k\cong E(\lambda_1)\otimes
P(v_2) \otimes P_{p-1}(b)\otimes \F_p\{a_1t^{dp}\,|\,0<d<p\}\,.
\end{equation*}
Thus $\ker(R-1)_*$ is a free $P(v_2)$-module, and the
exact sequence~\eqref{kku-eq:co-ker} splits. We have an isomorphism
of $P(v_2)$-modules
\begin{equation} 
\label{kku-eq:pre-tc}
\begin{split}
V(1)_*TC(ku_p)\cong\ & P(v_2)\otimes P_{p-1}(b) 
\otimes E(\partial, \lambda_1, a_1)
\\
&\oplus P(v_2)\otimes P_{p-1}(b)\otimes E(a_1)\otimes
\F_p\{\lambda_1t^d\,|\,0<d<p\}\\
&\oplus P(v_2)\otimes P_{p-1}(b)\otimes
E(\lambda_1)\otimes\F_p\{a_1t^{pd}\,|\,0<d<p \}
\end{split}
\end{equation}
in degrees larger than $2p-2$, where the summand 
\begin{equation*}
P(v_2)\otimes P_{p-1}(b) \otimes E(\lambda_1,
a_1)\otimes\F_p\{\partial\}
\end{equation*}  
is the group $\cok(R-1)_*\cong\Sigma^{-1}\ti A$.
\par
We now show that the relation
\begin{equation*}
b^{p-1}=-v_2
\end{equation*}
holds in $V(1)_*TC(ku_p)$. We recall from
Proposition~\ref{kku-prop:prebott} that the class $b'^{p-1}+v_2$ in
$V(1)_{2p^2-2}K(\Z,3)$ is annihilated by
$b'$. This class maps by $\trc_*\phi_*$ to the class
\begin{equation*}
b^{p-1}+v_2
\in V(1)_{2p^2-2}TC(ku_p)
\,,
\end{equation*}
which is therefore annihilated by $b$. 
Thus it suffices to show that zero is the only class in
$V(1)_{2p^2-2}TC(ku_p)$ 
that is annihilated by $b$.  
We consider the short exact
sequence
\begin{equation*}
  0\to \cok(R-1)_{2p^2-1}\to V(1)_{2p^2-2}TC(ku_p)\to
  \ker(R-1)_{2p^2-2}\to 0
\end{equation*}
given in~\eqref{kku-eq:co-ker} above. Here 
\begin{equation*}
\ker(R-1)_*\subset V(1)_*TF(ku_p)
\end{equation*}
inherits via $\Gamma_*$ the spectral sequence filtration  
of $V(1)_*\thh(ku_p)^{hS^1}$. 
By~\eqref{kku-eq:pre-tc}, 
this filtration gives the short exact sequence 
\begin{equation*}
  0\to \F_p\{{b^{p-2}\cdot \lambda_1 \cdot a_1t^p}\}\to
  \ker(R-1)_{2p^2-2}\to \F_p\{\overline{v_2}\}\to 0
\end{equation*} 
in dimension $2p^2-2$, 
while in dimension $2p^2+2p$ it gives
the short exact sequence 
\begin{equation*}
  0\to \F_p\{{v_2\cdot \lambda_1\cdot a_1t^p}\}\to
  \ker(R-1)_{2p^2+2p}\to \F_p\{\overline{b\cdot v_2}\}\to
  0\,.
\end{equation*} 
Here $\overline{v_2}$ and $\overline{b\cdot v_2}$ 
are represented by ${t\mu}$ and
and $b_1\cdot t\mu$ in $E^\infty(S^1)$, respectively. 
Multiplication with $b$ is compatible with the filtration,
and maps the former sequence to the latter one.
First, notice that the class $\overline{v_2}$ maps to a non-zero class in
$\F_p\{\overline{b\cdot v_2}\}$,
since $b\cdot \overline{v_2}$ is represented by
$b_1\cdot t\mu$ in $E^\infty(S^1)$. Next,
the relation $b^p=-bv_2$ in $\ker(R-1)_*$ implies
\begin{equation*}
b^{p}\cdot \lambda_1\cdot a_1t^p=-v_2\cdot b\cdot \lambda_1\cdot
a_1t^p\,,
\end{equation*}
which is non-zero by~\eqref{kku-eq:pre-tc}. A fortiori
$b^{p-1}\cdot \lambda_1\cdot a_1t^p\in 
\F_p\{v_2\cdot \lambda_1\cdot a_1t^p\} $ is not zero either.
Thus $\ker(R-1)_{2p^2-2}$ contains no
non-zero class annihilated by $b$, and we deduce that 
\begin{equation*}
  b^{p-1}+v_2\ \in\ 
  \partial(\cok(R-1)_{2p^2-1})\ =\ \F_p\{b^{p-2}\cdot
a_1\cdot\partial\}\,.
\end{equation*}
However the class $b^{p-2}\cdot a_1\cdot\partial$
is not annihilated by $b$, since by~\eqref{kku-eq:pre-tc} we
know that
$b^p\cdot a_1\cdot \partial =-v_2\cdot b\cdot a_1 \cdot\partial$
is non-zero.  This proves that $b^{p-1}+v_2$ must be zero.
\par
In particular $b$ is not a nilpotent class, and 
we have an isomorphism of $P(b)$-modules
\begin{equation*}
\begin{split}
V(1)_*TC(ku_p)\cong\ & P(b) \otimes E(\partial, \lambda_1, a_1)
\\
&\oplus P(b)\otimes E(a_1)\otimes
\F_p\{t^d\lambda_1\,|\,0<d<p\}\\
&\oplus P(b)\otimes
E(\lambda_1)\otimes\F_p\{a_1t^{pd}\,|\,0<d<p \}
\end{split}
\end{equation*}
in degrees larger than $2p-2$. This proves that our 
formula for $V(1)_*TC(ku_p)$ is correct
in dimensions greater than $2p-2$.
Let us define $M$ and $N$ as
\begin{equation*}
M=\kern -7pt\bigoplus_{-1\leq n\leq 2p-2}\kern -7ptV(1)_nTC(ku_p)\  \textup{ and }\
N=\kern -2pt\bigoplus_{n\geq 2p-1}\kern -2pt V(1)_nTC(ku_p)\,.
\end{equation*}
We just argued that $N$ is a free $P(b)$-module.
We know by~\eqref{kku-eq:TCS3} that there is an isomorphism 
\begin{equation*}
M\cong\F_p\{\partial,1,u^ia_0,\lambda_1t^d,\partial\lambda_1\,|\,
0\leq i\leq p-3,\ 1\leq d\leq p-1\}
\end{equation*}
of $\F_p$\,-modules.
This proves that the formula for $V(1)_*TC(ku_p)$ in Theorem~\ref{kku-thm:tcku} 
holds as an isomorphism of
$\F_p$-modules. It only remains to show 
that for any non-zero class $m\in M$, we
have $bm\neq0$ in $V(1)_*TC(ku_p)$. 
By comparison with $V(1)_*TC(\ell_p)$ or with
$V(1)_*\thh(ku_p)^{hS^1}$, we know that either $m\lambda_1$
or $mv_2$ is non-zero. These products lie in $N$ for
degree reasons, so are not $b$-torsion classes. Therefore $m$ is
not a $b$-torsion class either.
\end{proof}
\section{Algebraic K-theory}
\label{kku-sec:algK}
\begin{theorem}\label{kku-thm:kku}
There is an isomorphism of $P(b)$-modules
\begin{equation*}
\begin{split}
V(1)_*K(ku_p)\cong\ & P(b)\otimes E(\lambda_1,a_1)\oplus 
P(b)\otimes \F_p\{\partial\lambda_1,\partial b,\partial
a_1,\partial\lambda_1a_1\}\\
&\oplus  P(b)\otimes E(a_1)\otimes \F_p\{t^d\lambda_1\,|\,0<d<p \}\\
&\oplus  P(b)\otimes
E(\lambda_1)\otimes\F_p\{\sigma_n,\lambda_2t^{p^2-p}\,|\,1\leq
n\leq p-2 \}\\
&\oplus \F_p\{s\}\,,
\end{split}
\end{equation*}
with $b^{p-1}=-v_2$. The degree of the generators is given by 
$|\partial|=-1$, $|\lambda_1|=2p-1$, $|a_1|=2p+3$,
$|\sigma_n|=2n+1$,
$|t|=-2$, $|\lambda_2|=2p^2-1$ and $|s|=2p-3$.
The classes $1$, $\sigma_n$, $\lambda_1$, $b$ and $a_1$ map under
the trace to $1$, $u^{n-1}a_0$, $\lambda_1$, $b_1$ and $a_1$ in
$V(1)_*\thh(ku_p)$, respectively,
and the other given $P(b)$-module generators map to zero. 
\end{theorem}
\begin{proof}
There is a cofibre sequence of spectra~\cite{Rkku-HM1}
\begin{equation*}
K(ku_p)_p\to TC(ku_p)\to \Sigma^{-1}H\Z_p\to \Sigma
K(ku_p)_p\,.
\end{equation*}
We have an isomorphism $V(1)_*\Sigma^{-1}H\Z_p\cong \F_p\{\partial,\epsilon\}$ 
with a primary $v_1$ Bockstein $\beta_{1,1}(\epsilon)=\partial$. 
Here $\partial$ is the image of the class
$\partial\in V(1)_{-1}TC(ku_p)$, while $\epsilon$ maps by the
connecting homomorphism to a class $s\in V(1)_{2p-3}K(ku_p)$. 
These facts, together with Theorem~\ref{kku-thm:tcku},
allow us to establish our formula for $V(1)_*K(ku_p)$. 
The statement on the trace follows from the definition of
the given $P(b)$-module generators. 
\end{proof}
\par
The following corollary is a restatement of
Proposition~\ref{kku-prop:desc} part (b) of the
introduction.
\begin{corollary}\label{kku-cor:lvsku}
There is a short exact sequence of $P(b)$-modules
\begin{equation*}
0\to K\to
P(b)\otimes_{P(v_2)}V(1)_*K(\ell_p)\stackrel{\mu}{\lr} K(ku_p) \to Q \to 0
\end{equation*}
where $K$ and $Q$ are finite (and hence torsion)
$P(b)$-modules given by
\begin{equation*}
  \begin{split}
    K=\ &\F_p\{b^{k}a\,|\, 1\leq k\leq p-2\},\ \textup{and}\\
    Q=\ &P_{p-2}(b)\otimes\F_p\{ \partial b, \partial a_1, 
    a_1, \partial\lambda_1a_1, \lambda_1a_1 \}\\
    &\oplus P_{p-2}(b)\otimes\F_p\{a_1\lambda_1t^d\,|\,0<d<p\}\\
    &\oplus E(\lambda_1)\otimes\F_p\{\sigma_nb^{i_n}\,|\,1\leq
    n\leq p-2,\ 0\leq i_n\leq p-2-n\}\,.
  \end{split}
\end{equation*}
Here $a\in V(1)_{2p-3}K(\ell_p)$ is the class annihilated
by $v_2$ and mapping to $s$. 
In particular we have an isomorphism $P(b,b^{-1})$-algebras  
\begin{equation*}
P(b,b^{-1})\otimes_{P(v_2)}V(1)_*K(\ell_p)\cong
V(1)_*K(ku_p)[b^{-1}]\,. 
\end{equation*}
\end{corollary}
\begin{proof}
 This follows from the formulas for $V(1)_*K(\ell_p)$
 and for $V(1)_*K(ku_p)$ given in~\cite{Rkku-AR02}*{Theorem~9.1} 
and Theorem~\ref{kku-thm:kku}, and the fact that
 $V(1)_*K(\ell_p)$ includes as the summand of
 $\delta$-weight zero in $V(1)_*K(ku_p)$,
 see~\cite{Rkku-Au05}*{Theorem~10.2}. Notice that for $1\leq
 d\leq p-2$ the class
 $\lambda_2t^{dp}\in V(1)_{2p^2-pd-1}K(\ell_p)$ maps to
 $\sigma_db^{p-1-d}$, up to a non-zero scalar multiple. 
\end{proof}
\par
Blumberg and Mandell~\cite{Rkku-BlM} have proved a conjecture of John Rognes
that there is a localization cofibre sequence 
\begin{equation*}
K(\Z_p) \stackrel{\tau}{\lr}  K(ku_p) \stackrel{j}{\lr}
K(KU_p)\to \Sigma  K(\Z_p)\,,
\end{equation*}
relating the algebraic
$K$-theory of $ku_p$, of its localization $KU_p=ku_p[u^{-1}]$ (i.e.
periodic $K$-theory), and of its mod $(u)$ reduction $H\Z_p$.
The $V(1)$-homotopy of $K(\Z_p)$ and $K(ku_p)$ is known, but we need to
compute also the transfer map $\tau_*$ and solve a $P(b)$-module
extension if we seek a decent description of $V(1)_*K(KU_p)$. 
Let us therefore assume that this localization sequence maps
via trace maps to a
corresponding localization sequence in topological Hochschild
homology, building a homotopy commutative diagram of horizontal fibre
sequences
\begin{equation}\label{kku-eq:loc-seq}
\xymatrix@R=6ex@C=4ex{
K(\Z_p)\ar[d]^-{\mathrm{tr}}\ar[r]^{\tau}&
K(ku_p)\ar[d]^-{\mathrm{tr}}\ar[r]^-j &
K(KU_p)\ar[d]^-{\mathrm{tr}}\ar[r] &
\Sigma K(\Z_p)\ar[d]^-{\Sigma\mathrm{tr}}\\
\thh(\Z_p)\ar[r]^{\tau}&
\thh(ku_p)\ar[r]^-j &
 \thh(ku_p|KU_p)\ar[r] &
\Sigma \thh(\Z_p)  \,,
}
\end{equation}
as conjectured by Lars Hesselholt, compare with
Remark~\ref{kku-rem:HMlog} below. 
The $V(1)$-homo\-topy of the bottom line was described
in~\cite{Rkku-Au05}*{\S10}.
The $V(1)$-homotopy groups of
$K(\Z_p)$ are given by an isomorphism~\cite{Rkku-HM1}
\begin{equation*}
V(1)_*K(\Z_p)\cong E(\lambda_1)\oplus \F_p\{\partial v_1,
\partial\lambda_1\}
\oplus \F_p\{ \lambda_1t^{d}\,|\,0<d<p \}\,.
\end{equation*}
The class $\partial v_1$ maps to $s$ in $V(1)_*K(ku_p)$ via
$\tau_*$.
The class $1\in V(1)_0K(\Z_p)$ is in the kernel of ${\tau}_*$, because it is
$v_2$-torsion and there is no torsion class in $V(1)_0K(ku_p)$.
Let $d\in V(1)_1K(KU_p)$ be the class mapping 
to $1\in V(1)_0K(\Z_p)$ via the connecting homomorphism.
Presumably $d$ corresponds to the added unit or the self-equivalence 
\begin{equation*}
KU_p\stackrel{u}{\lr} \Sigma^{-2}KU_p \stackrel{\simeq}{\lr}KU_p\,,
\end{equation*}
where $u$ denotes multiplication by the Bott class, and the second
map is the Bott equivalence.
The class $d$ maps in $V(1)_1\thh(ku_p|KU_p)$ to a class
with the same name.
In~\cite{Rkku-Au05}*{\S10} we
establish an (additive) isomorphism
\begin{equation}\label{kku-eq:log-thh}
V(1)_*\thh(ku_p|KU_p)\cong P_{p-1}(u)\otimes E(d,\lambda_1)\otimes
P(\mu_1)\,.
\end{equation}
If this is an isomorphism of algebras, then 
$j_*(b_1)d=j_*(a_1)$ holds in $V(1)_*\thh(ku_p|KU_p)$,  
and lifts to the relation
$j_*(b)d=j_*(a_1)$ in $V(1)_*K(KU_p)$. By inspection this determines the structure of
$V(1)_*K(KU_p)$ as a $P(b)$-module.
\begin{theorem}\label{kku-thm:periodic}
Under the hypothesis that there
exists a commutative diagram of localization
sequences~\textup{\eqref{kku-eq:loc-seq}},
and that the isomorphism~\textup{\eqref{kku-eq:log-thh}} 
is one of algebras, we have an isomorphism of $P(b)$-modules
\begin{equation*}
\begin{split}
V(1)_*K(KU_p)\cong\ & P(b)\otimes E(\lambda_1,d)\oplus 
P(b)\otimes \F_p\{\partial\lambda_1,\partial b,\partial
a_1,\partial\lambda_1d\}\\
&\oplus  P(b)\otimes E(d)\otimes \F_p\{t^d\lambda_1\,|\,0<d<p \}\\
&\oplus  P(b)\otimes
E(\lambda_1)\otimes\F_p\{\sigma_n,\lambda_2t^{p^2-p}\,|\,1\leq
n\leq p-2
\}\,.
\end{split}
\end{equation*} 
The class $d$ has degree $1$, and the other classes have the
degree given in Theorem~\textup{\ref{kku-thm:kku}}.
\end{theorem}
\begin{remark}\label{kku-rem:HMlog}
Consider a complete discrete valuation field $K$ of characteristic
zero with perfect residue field $k$ of characteristic
$p\geq3$, and let $A$ be its valuation ring. 
Hesselholt and Madsen~\cite{Rkku-HM2} compute the 
$V(0)$-homotopy of $K(A)$ and $K(K)$ by means of the
cyclotomic trace. They introduce a relative version
of topological cyclic homology, denoted $TC(A|K)$, that sits in 
a localization cofibre sequence
\begin{equation*}
  TC(k)\to TC(A)\to TC(A|K)\to \Sigma TC(k)\,.
\end{equation*}
The computation of 
$V(0)_*TC(A|K)$ is achieved by using the rich algebraic
structure on the $V(0)$-homotopy groups of the tower 
$TR^\bullet(A|K)$, and described in terms of the de\,Rham-Witt complex with log poles
\begin{equation*}
W_\bullet\omega^*(A,A\cap K^\times)\,,
\end{equation*}
see~\cite{Rkku-HM2}*{Th.\,C}.  
Then $V(0)_*TC(A)$ can be evaluated by means of the
localization sequence. This approach has, in particular, the advantage of
avoiding a computation of $V(0)_*TR^\bullet(A)$, which seems quite intractable.
\par
Continuing the discussion in~\cite{Rkku-Au05}*{\S10} on a relative trace for $ku_p$, 
and following Lars Hesselholt, 
one could speculate on the existence of a relative term $TC(ku_p|KU_p)$ 
fitting in a localization sequence
\begin{equation*}
  TC(H\Z_p)\to TC(ku_p)\to TC(ku_p|KU_p)\to \Sigma
  TC(H\Z_p)\,,
\end{equation*}
through which the trace of diagram~\eqref{kku-eq:loc-seq} factorizes.
By analogy with the case of complete discrete valuation
fields, we expect that a computation  
of $V(1)_*TR^n(ku_p|KU_p)$ 
should be easier to handle than the computation
of $V(1)_*TR^n(ku_p)$ presented in this paper.
In fact, the advantage of such an approach is already apparent when comparing 
\begin{equation*}
V(1)_*TR^1(ku_p|KU_p)=V(1)_*\thh(ku_p|KU_p)
\end{equation*} 
in~\eqref{kku-eq:log-thh} with $V(1)_*\thh(ku_p)$
in~\eqref{kku-eq:thh}, and is also confirmed by
partial, hypothetical computations of
$V(1)_*TR^n(\ell_p|L_p)$ and
$V(1)_*TR^n(ku_p|KU_p)$ by Lars Hesselholt (private communication) 
and the author.  
\end{remark}
\begin{acknowledgements}
This paper is part of my Habilitation
thesis written at the University of Bonn. I thank Stefan Schwede, Carl-Friedrich
B{\"o}digheimer, G{\'e}rald Gaudens and my other colleagues
in Bonn for their friendly support. 
I thank Birgit Richter, Bj{\o}rn Dundas, John Greenlees, Lars Hesselholt,
Christian Schlichtkrull and Neil Strickland for interesting conversations related to
this project. 
This paper builds on the results of~\cite{Rkku-AR02}, and uses many 
techniques and ideas that I learned from John Rognes. I am very grateful 
to him for his help and his generosity. Finally, I
thank the referee for his many useful suggestions.
\end{acknowledgements}

\end{document}